\documentclass[oneside,english]{amsart}
\usepackage[T1]{fontenc}
\usepackage[latin9]{inputenc}
\usepackage{amsthm}
\usepackage{graphicx}
\usepackage{hyperref}

\makeatletter
\numberwithin{equation}{section}
\numberwithin{figure}{section}

\makeatother

\usepackage{babel}

\newtheorem{theorem}{\bf Theorem}[section]
\newtheorem{lemma}{\bf Lemma}[section]

\newtheorem{assumption}{\bf Assumption}

\begin{document}

\title[Asymptotic solution to convolution integral equations]{Asymptotic solution to convolution integral equations on large and small intervals}

\author{
Dmitry Ponomarev$^{1,2}$ \email{dmitry.ponomarev@asc.tuwien.ac.at}}

\address{$^{1}$TU Wien, Wiedner Hauptstrasse 8, 1040 Wien, Austria\\ $^{2}$St. Petersburg Department of Steklov Mathematical Institute RAS, Fontanka 27, 191023 St. Petersburg, Russia}



\begin{abstract}
We consider convolution integral equations on a finite interval with a real-valued kernel of even parity, a problem equivalent to finding a Wiener-Hopf factorisation of a notoriously difficult class of $2\times 2$ matrices. The kernel function is assumed to be sufficiently smooth and decaying for large values of the argument. Without loss of generality, we focus on a homogeneous equation and we propose methods to construct explicit asymptotic solutions when the interval size is large and small. The large interval method is based on a reduction of the original equation to an integro-differential equation on a half-line that can be asymptotically solved in a closed form. This provides an alternative to other asymptotic techniques that rely on fast (typically exponential) decay of the kernel function at infinity which is not assumed here. We also consider the problem on a small interval and show that finding its asymptotic solution can be reduced to solving an ODE. In particular, approximate solutions could be constructed in terms of readily available special functions (prolate spheroidal harmonics). Numerical illustrations of the obtained results are provided and further extensions of both methods are discussed. 
\end{abstract}


\maketitle

\section{Introduction}
\label{sec:intro}

One-dimensional convolution integral equations on finite intervals are ubiquitous in applied and theoretical physical contexts and appear in a number of mathematical and engineering areas \cite{Gr64, vTr69, Ba73, FaLa20, LePo17}. Namely, we consider the Fredholm integral equation
\begin{equation}
\int_{A}K\left(x-t\right)f\left(t\right)dt=\lambda f\left(x\right)+g\left(x\right),\hspace{1em}x\in A,\label{eq:main_eq_nonhom}
\end{equation}
where $A:=\left(-a,a\right)$, $a>0$, is a finite interval, and $K\left(x\right)$, $g\left(x\right)$ are given functions, $\lambda$ is a constant, and $f\left(x\right)$ is the solution to be found.

We note immediately that due to the convolution structure of the kernel, the equation is translationally invariant, so the convenient choice of the symmetric interval does not restrict generality of the problem.   

If $A$ was the whole real line, the solution could be found simply by taking Fourier transform and solving an algebraic equation in the Fourier domain. If $A$ was a half-line, the equation would be of a Wiener-Hopf type and hence still solvable in an explicit form \cite{No88, GaCh78} though the solution method is more involved. In the present case, when $A$ is a finite interval, the Wiener-Hopf theory could still be useful and, from that viewpoint, solving (\ref{eq:main_eq_nonhom}) is tantamount to factorisation of a triangular $2\times 2$ matrix with non-factorisable (exponential) diagonal elements. Contrary to the scalar Wiener-Hopf factorisation arising in half-line problems, explicit constructions of matrix factorisations is generally available only for certain clases of matrices \cite{RoMi16} and dealing with the class of matrices with diagonal exponential factors is not an easy business, even on a qualitative level \cite{BoSp13}. This difficulty reflects the fact that, in general, convolution integral equations on a finite interval cannot be solved in a closed form and, therefore, one has to resort to approximation methods such as numerical discretisation, iterative solution or asymptotic methods. 

Considering problem (\ref{eq:main_eq_nonhom}) with a real-valued non-singular symmetric kernel $K(x-t)$, we focus on the homogeneous equation
\begin{equation}
\int_{A}K\left(x-t\right)f\left(t\right)dt=\lambda f\left(x\right),\hspace{1em}x\in A.\label{eq:main_eq}
\end{equation}

We thus deal with a problem of spectral decomposition of a convolution integral operator on an interval. In view of completeness of the set of eigenfunctions in the range of the correponding integral operator, problem \eqref{eq:main_eq} is the most general in a sense that its solution also allows treating inhomogeneous Fredholm equations of both first and second kinds. 
 
On the other hand, compared to its inhomogeneous counterpart (\ref{eq:main_eq_nonhom}), the problem (\ref{eq:main_eq}) is more delicate as it requires more sophisticated analysis due to the existence of generally infinite number of solutions (eigenfunctions) and involvement of an unknown spectral parameter (eigenvalue).

The problem \eqref{eq:main_eq_nonhom} is certainly not new and has been previously addressed in numerous works. We mention here only papers which are constructive meaning that they describe a way of computing solutions rather than studying their qualitative features (such as decay to zero of eigenvalues and their bounds, cf. \cite{BeLa52, GrSz56, Wi64}).
In \cite{La58}, assuming even parity and exponential decay of the kernel function at infinity, the author employs Wiener-Hopf factorisation in a strip to reformulate the problem as a pair of integral equations in Fourier domain which are further solved approximately for the large interval case.
Hutson \cite{Hu64} builds up on the work of Latter, simplifying and justifying the approximation results rigorously under an additional assumption that the Fourier transform on the kernel function is positive and decreases strictly monotonically with the distance from the origin.
In a duplet of papers \cite{LeMu1, LeMu2} it is assumed that the kernel function admits certain representation in the form of a Laplace integral and the problem is treated through reduction to integral equations for the so-called Chandrasekhar's $X$- and $Y$-functions (their analogs are also known as Ambarzumian functions). Even though these auxiliary equations cannot be solved explicitly, the authors claim an implication that solutions of the homogeneous version of the original integral equation are sine and cosine functions with some frequencies that are not available in a closed form.
Revisiting an earlier result of Widom \cite{Wi61} for the large interval approximation of solutions of \eqref{eq:main_eq}, van Trigt in \cite{vTr73} discusses reformulations into a canonical integral equation in Fourier domain whose kernel depends only on the behavior of Fourier transform of the original kernel function near the origin.
Making suitable assumptions on the kernel function, the method of formal asymptotic expansion has been used in \cite{Ga89} to arrive at a series asymptotic to the solution when the interval is large.
In \cite{OlGa89, KnKe91}, the formal asymptotic expansion technique (pioneered for integral equations by Carrier \cite{Ca55}) is combined with a boundary layer solution obtained in vicinity of the endpoints upon rescaling. The solution is expressed in terms of moments of the kernel function on the entire line which are assummed to be finite and this hence entails the decay assumption at infinity such that at least $\int_{\mathbb{R}}x^2 K\left(x\right)dx<\infty$. For the case of the small interval, the authors of \cite{KnKe91} also construct approximations of the solution in a polynomial form providing the kernel function is differentiable sufficiently many times. 
We also mention few works \cite{GoFe74, Sa13, Po91, Ga84} which, without any assumptions on the size of interval, provide representations of the solution of \eqref{eq:main_eq_nonhom} for a general function $g$ in terms of solutions for its specific choices such as the kernel function itself or an exponential function, i.e. essentially reducing the problem to an equation for the resolvent.
Another approach geared for inhomogeneous equations is based on treating the problem on an interval of variable size using the so-called invariant imbedding method \cite{KaKa68} that, under assumptions on the kernel function similar to those in \cite{LeMu1}, allows obtaining solution by means of numerical integration of a system of nonlinear integro-differential equations.
  
In the present work, we propose a method that is able to overcome the requirement of exponential decay of the kernel function at infinity providing an explicit large interval approximation to eigenvalues and eigenfunctions. We assume only mild smoothness and algebraic decay of the kernel function in addition to another typical assumption (see \cite{Hu64,vTr73}) of positivity and monotone decay of the Fourier transform of the kernel on the positive half-line. In contrast with other works, assumptions are thus made only on behavior of the Fourier transform on the kernel function on the real line and not in the entire complex domain and/or in a strip containing the real line. This improvement is not just an academic curiosity but strongly motivated by kernel functions such as the families of functions $\frac{1}{\left(x-t\right)^2+h^2}$, $\frac{1}{\left(\left(x-t\right)^2+h^2\right)^{3/2}}$ defined for a parameter $h>0$ and arising in different practical and theoretical contexts (see e.g. \cite{FaLa20, Ki11, LePo17}).
When the interval $A$ is small, we provide an alternative to the conventional polynomial approximation of the kernel function. We present an approach that hinges on the approximation of the original operator by another compact operator of a non-finite rank thus preserving a structural property of the problem.

Employing general properties of the solution given in Section \ref{sec:asss_and_props}, our large interval approach described in Section \ref{sec:large_A} reduces the original equation to a problem for determining continuation of the solution onto a half-line. Such a continuation satisfies an integro-differential equation on a half-line with two kernel functions depending on sum and difference of the arguments. For the large interval, the dominant dependence is that of the difference type, a fact that allows solving this approximate half-line problem explicitly reformulating it as a Riemann-Hilbert problem. The solution inside the interval is recovered from the half-line continuation that has thus been approximately found. Matching of the the original solution and its half-line continuation results in a pair of characteristic equations that allows finding eigenvalues corresponding to even and odd solutions (eigenfunctions). Approximate solutions of the original problem are furnished as the sum of a trigonometric function (sine or cosine) perturbed by an integral term given explicitly. The result is formulated in the form of a representation theorem of the exact solution provided implicitly whereas the approximation issues are discussed afterwards and suggest an explicit form of the approximate solution.
In Section \ref{sec:small_A}, we treat the case of a small interval size by advocating a novel approach which consists in approximating the original equation by an integral equation whose kernel admits a commuting differential operator and hence reduces a problem to integrating an ODE. In particular, we point out that the original equation can be solved in terms of prolate spheroidal wave functions which can be effectively computed. 
Section \ref{sec:num_ill} illustrates the proposed asymptotic methods on a particular equation and compares eigenfunctions and eigenvalues with reference solutions obtained by a direct numerical discretisation.
Finally, in Section \ref{sec:concl}, we summarise the obtained results and discuss their possible extensions.

\section{Assumptions and solution properties}
\label{sec:asss_and_props}

Let us set some notation that we shall use throughout the paper. 
We employ notation $\chi_S\left(x\right)$ for the characteristic function of a set $S$. We use $\bar{S}$ to mean the closure of the set $S$. We denote $\mathbb{R}_{+}$, $\mathbb{R}_{-}$ positive and negative real half-lines, and $\mathbb{Z}$, $\mathbb{N}_{0}$, $\mathbb{N}_{+}$ the sets of all integer numbers and those starting from $0$ and $1$, respectively.  We resort to the Lebesgue measure notation and write $\left|A\right|$ meaning the size of the interval $A$. The standard notation for the derivatives is used meaning that $f^{\prime}$, $f^{\prime\prime}$, $f^{\left(j\right)}$ are the first, the second and the $j$-th derivative of a function $f$, respectively. We frequently use the "big-O" and "small-O" Landau notations, that is, given functions $f$ and positive-valued $g$, we write $f\left(x\right)=\mathcal{O}\left(g\left(x\right)\right)$ for $x\rightarrow\infty$ (or equivalently $x\gg 1$) meaning that for arbitrary large $x$ we have the bound $\left|f\left(x\right)\right|\leq C g\left(x\right)$ with some constant $C>0$, and we write $f\left(x\right)=o\left(g\left(x\right)\right)$ if $\displaystyle{\lim_{x\rightarrow\infty}\left|f\left(x\right)\right|/g\left(x\right)}=0$ for $x\rightarrow\infty$. We abbreviate $f\left(x\pm 0^{+}\right)=\displaystyle{\lim_{\epsilon\rightarrow 0}}f\left(x\pm\epsilon\right)$. To stress the positivity of $\epsilon$, we sometimes write $\epsilon\searrow 0$ or $\epsilon\rightarrow 0^{+}$ instead of simply $\epsilon\rightarrow 0$. We use $\text{sgn }x$ notation for the sign function which takes value $1$ if $x$ is positive and value $-1$ otherwise. We write $\left\lfloor X\right\rfloor$ to denote the integer part of a number $X$ and use $\overline{X}$ for the complex conjugate of $X$. We use ":=" sign to introduce new quantities. We write $f\left(x\right)\equiv C$ to stress that a function $f\left(x\right)$ is identically constant (and equal to $C$). The notation "$\text{p.v.} \int$" stands for the Cauchy principal value integral. Our Fourier transform convention is $\hat{f}\left(k\right):=\int_\mathbb{R}f\left(x\right)e^{2\pi ikx}dx$. We denote Hilbert transform of $f$ as $\mathcal{H}\left[f\right]\left(x\right):=\frac{1}{\pi}\text{p.v.}\int_{\mathbb{R}}\frac{f\left(\tau\right)}{x-\tau}d\tau$ and introduce the following operators $P_{\pm}\left[f\right]:=\frac{1}{2}f\pm \frac{i}{2}\mathcal{H}\left[f\right]$. Given $0<\gamma\leq 1$, $j\in\mathbb{N}_{0}$, $1\leq p\leq\infty$, we denote $C^{0,\gamma}$, $C^j$, $L^p$ the spaces of H\"{o}lder continuous functions, functions whose $j$-th derivative is continuous and the Lebesgue space for $1\leq p\leq\infty$, respectively. Note that when $\gamma=1$, $C^{0,\gamma}=C^{0,1}$ is the space of Lipschitz continuous functions. 

Let us denote $\mathfrak{K}$ the convolution integral operator in the left-hand side of (\ref{eq:main_eq}).

We need the following assumptions on the kernel function of (\ref{eq:main_eq}) some of which have already been mentioned above.

\begin{assumption}\label{ass:K_par}
$K$ is an even real-valued function: $K\left(x\right)=K\left(-x\right)\in\mathbb{R}$ for $x\in\mathbb{R}$.
\end{assumption}
\begin{assumption}\label{ass:K_reg}
$K\in C^{2}\left(\mathbb{R}\right)$, $K^{\left(j\right)}\left(x\right)=\mathcal{O}\left(1/x^\alpha\right)$ as $\left|x\right|\rightarrow\infty$ for some $\alpha>1$ and any $j\in\left\{0,1,2\right\}$,  $\hat{K}\in C^1\left(\mathbb{R}\backslash\left\{0\right\}\right)$ and, moreover, there exist finite $\hat{K}^\prime\left(0^-\right)$, $\hat{K}^\prime\left(0^+\right)$.
\end{assumption}
\begin{assumption}\label{ass:pos_def}
$\mathfrak{K}$ is positive definite, that is for any $f\in L^2\left(A\right)$ we have
\begin{equation}
\left<\mathfrak{K}f,f\right>_{L^2\left(A\right)}:=\int_\mathbb{R}\int_A K\left(x-t\right)f\left(t\right)\overline{f\left(x\right)}dtdx>0.
\end{equation}
\end{assumption}
\begin{assumption}\label{ass:FT_monot_decay}
$\hat{K}\left(k\right): \mathbb{R}\rightarrow \mathbb{R}_{+}$ decays strictly monotonically with $\left|k\right|$.
\end{assumption}

Some comments about the made assumptions are now in place.
There is a certain redundancy in the conditions of Assumption \ref{ass:K_reg} as they deal with both regularity and decay of a function and its Fourier transform and those are not independent.  
Note that the regularity and the decay condition for $K$ imply that $\hat{K}\in C\left(\mathbb{R}\right)$ and hence Assumption \ref{ass:K_reg} means that $\hat{K}\in C^{0,1}\left(\mathbb{R}\right)$ with the origin $x=0$ being the only point of non-differentiability. 
It goes without saying that the positive definiteness in Assumption \ref{ass:pos_def} can be replaced by the negative definiteness while simultaneously changing the monotone decay to the monotone growth condition in Assumption \ref{ass:FT_monot_decay}.
Assumption \ref{ass:K_reg} will only be relevant for Section \ref{sec:large_A}, in Section \ref{sec:small_A} regularity assumptions will be discussed separately. 

Elementary properties of the solutions of \eqref{eq:main_eq} are described in
\begin{lemma}\label{lem:props}
Suppose that $K$ satisfies Assumptions \ref{ass:K_par}-\ref{ass:FT_monot_decay}.
 Then, there exists a set $\left(\lambda_{n}\right)_{n\geq 1}$ such that for each $\lambda_n$ there is only one $f_n\in C^2\left(\bar{A}\right)$ which satisfies (\ref{eq:main_eq}) and is either an odd or an even function. Moreover, $0<\ldots<\lambda_{n+1}<\lambda_n<\ldots<\lambda_1<\hat{K}\left(0\right)$, the set $\left(f_n\right)_{n\geq 1}$ forms an orthogonal basis in the range of $\mathfrak{K}$, and $\left|f_n\left(a\right)\right|=\left|f_n\left(-a\right)\right|\neq 0$. Without loss of generality, each $f_n$ can be chosen to be real-valued and normalised such that $\left\Vert f_n\right\Vert_{L^2\left(A\right)}=1$.
\end{lemma}

\begin{proof}
Since $A$ is finite and $K$ is regular (Assumption \ref{ass:K_reg}), the operator $\mathfrak{K}$: $L^{2}\left(A\right)\rightarrow L^{2}\left(A\right)$ is compact as a Hilbert-Schmidt operator \cite[Thm 3.2.7]{Ha95}. By Assumption \ref{ass:K_par}, $\mathfrak{K}$ is also self-adjoint. 
By the spectral theorem for compact self-adjoint operators \cite[Sect 6.11]{NaSe82}, the spectrum of $\mathfrak{K}$ is a non-empty purely discrete set of real eigenvalues $\left(\lambda_{n}\right)_{n\geq1}$ with zero being the only possible accumulation point. For each $\lambda_n$ there are potentially several but finitely many (up to some $M_n\in\mathbb{N}_{+}$) eigenfunctions and $\left(f_{n,m}\right)_{n\geq 1,1\leq m\leq M_n}\subset L^2\left(A\right)$ is a complete set in the range of $\mathfrak{K}$ such that eigenfunctions corresponding to different eigenvalues are orthogonal.

Since $\lambda\neq 0$, the smoothness of each eigenfunction follows from that of $K$ due to equation \eqref{eq:main_eq}.

The even parity and real-valuedness (Assumption \ref{ass:K_reg}) imply that $\hat K\left(k\right)$ is real-valued and even for $k\in\mathbb{R}$. Moreover, by Assumptions \ref{ass:pos_def}-\ref{ass:FT_monot_decay}, $\hat{K}\left(k\right)$ is positive for $k\in\mathbb{R}$ attaining its maximal value at $k=0$. Therefore, integrating both sides of \eqref{eq:main_eq} over $A$ against $\overline{f}$ and using the convolution theorem and isometry of Fourier transform (the Parseval's identity for inner products and norms), we obtain
\[
\lambda \left\Vert f\right\Vert_{L^2\left(A\right)}^2=\int_\mathbb{R}\left(\mathfrak{K}\chi_A f\right)\left(t\right)\overline{\chi_A\left(t\right)f\left(t\right)}dt=\int_\mathbb{R} \hat{K}\left(k\right)\left|\widehat{\chi_Af}\left(k\right)\right|^2 dk \in \left(0, \hat{K}\left(0\right)\left\Vert f\right\Vert_{L^2\left(A\right)}^2 \right)
\]
which implies that $0<\lambda_{n}<\hat{K}\left(0\right)$, $n\geq 1$.

The simplicity of the spectrum is not immediate but can be shown to hold under Assumptions \ref{ass:pos_def}-\ref{ass:FT_monot_decay}. The fact that for each eigenvalue $\lambda\in\left(\lambda_n\right)_{n=1}$ there is only one eigenfunction $f\left(x\right)$ was previously stated for particular kernels in \cite{BaLePo19, GoPa75, SlPo61}, however, the proof strategies extend to a general class of kernel functions. We provide such a proof here as it motivates Assumptions \ref{ass:pos_def}-\ref{ass:FT_monot_decay} and we could not find it in literature (though it would be naive to assume that this result is new given a vast amount of works on the general theory of Wiener-Hopf/Toeplitz operators).

The first part of the proof of the simplicity of the spectrum is to show that no eigenfunction $f\left(x\right)$ can have boundary values $f\left(-a\right)$, $f\left(a\right)$ vanishing simultaneously. This non-vanishing property can be deduced from $\hat{K}^\prime\left(k\right)\leq 0$, $k\geq 0$, and $\lambda>0$  using the following identity 
\begin{equation}
-\int_{0}^{\infty}k\hat{K}^{\prime}\left(k\right)\left[\left|\hat{f}_A\left(k\right)\right|^{2}+\left|\hat{f}_A\left(-k\right)\right|^{2}\right]dk=\lambda a\left[\left|f\left(a\right)\right|^{2}+\left|f\left(-a\right)\right|^{2}\right],\label{eq:energy}
\end{equation}
where $\hat{f}_A\left(k\right):=\int_A e^{2\pi i kx}f\left(x\right)dx.$

Relation \eqref{eq:energy} is reminiscent of an energy identity for differential equations in that it controls a solution in the interior of the domain in terms of the solution boundary values. 
It can be obtained as follows. First, recalling that $f$ and $K$ are smooth, we differentiate both sides of \eqref{eq:main_eq}, multiply by characteristic function $\chi_A$, take Fourier transform and use the convolution theorem to obtain
\begin{equation}
-2\pi i\int_{\mathbb{R}}\mathcal{S}_a\left(k,\tilde{k}\right)\tilde{k}\hat{K}\left(\tilde{k}\right)\hat{f}_A\left(\tilde{k}\right)d\tilde{k}=\lambda\int_{A}f^{\prime}\left(x\right)e^{2\pi ikx}dx,\hspace{1em} k\in\mathbb{R},\label{eq:energy_prelim1}
\end{equation}
where
\begin{equation}
{\mathcal S}_{a}\left(k,\tilde{k}\right):=\hat{\chi}_A\left(k-\tilde{k}\right)=\frac{\sin\left(2\pi a\left(k-\tilde{k}\right)\right)}{\pi\left(k-\tilde{k}\right)},\hspace{1em}k,\tilde{k}\in\mathbb{R}.\label{eq:RK_PW}
\end{equation}
Note that \eqref{eq:RK_PW} is a reproducing kernel for the Paley-Wiener space of exponent $2\pi a$ which is the space of Fourier transforms of functions supported on $\left(-a,a\right)$, i.e. a subspace of $L^2\left(\mathbb{R}\right)$ that consists of entire functions which are square-integrable over any horisontal line in $\mathbb{C}$ and obey the global bound $C e^{2\pi a\left|k\right|}$, $k\in\mathbb{C}$, for some constant $C>0$. If $\hat{f}_A$ belongs to this space, then so does $\overline{\hat{f}_A}$ and $\overline{\hat{f}^{\prime}_A}=-2\pi i \widehat{\chi_A x \overline{f}}$, and hence we have the identity
\begin{equation}
\int_{\mathbb{R}}{\mathcal S}_{a}\left(k,\tilde{k}\right)\overline{\hat{f}^{\prime}_A\left(\tilde{k}\right)}d\tilde{k}=\overline{\hat{f}^{\prime}_A\left(k\right)},\hspace{1em} k\in\mathbb{R}\label{eq:RKHS_ident}.
\end{equation}
Integrating both sides of \eqref{eq:energy_prelim1} against $\overline{\hat{f^{\prime}_A}}$ on $\mathbb{R}$ and employing \eqref{eq:RKHS_ident} yields
\begin{equation}
-2\pi i \int_{\mathbb{R}}k\hat{K}\left(k\right)\hat{f}_A\left(k\right)\overline{\hat{f}^{\prime}_A\left(k\right)}dk=\lambda\int_\mathbb{R}\overline{\hat{f}^{\prime}_A\left(k\right)}\int_A f^{\prime}_A\left(x\right)e^{2\pi ikx}dxdk,\label{eq:energy_prelim2}
\end{equation}
where in the left-hand side we employed Fubini's theorem to exchange the integration order. Indeed, since $\hat{f}_A\left(k\right)=\mathcal{O}\left(1/k\right)$ (as can be seen by integration by parts of the Fourier integral), it follows from smoothness of $\hat{f}_A$ that $k\hat{f}_A\left(k\right)<C_0$ for any $k\in\mathbb{R}$ and some constant $C_0>0$, and therefore, by the Cauchy-Schwarz inequality, we have 
\[
\int_\mathbb{R}\int_\mathbb{R}\left|\hat{f}_A^\prime\left(k\right)\right|\left|\frac{\sin\left(2\pi a\left(k-\tilde{k}\right)\right)}{\pi\left(k-\tilde{k}\right)}\right|dk\left|\tilde{k}\right|\left|\hat{f}_A\left(\tilde{k}\right)\right|\left|\hat{K}\left(\tilde{k}\right)\right|d\tilde{k}
\]\[
\leq C_0\left\Vert\hat{f}_A^\prime\right\Vert_{L^2\left(\mathbb{R}\right)}\left\Vert K\right\Vert_{L^1\left(\mathbb{R}\right)}\int_\mathbb{R}\left(\frac{\sin\left(2\pi a k\right)}{\pi k}\right)^2dk<\infty,
\]
which justifies the Fubini-Tonelli argument for the performed integral swap in \eqref{eq:energy_prelim2}.
After adding the complex conjugate and using the Parseval's identity for inner products in the right-hand side, equation \eqref{eq:energy_prelim2} transforms into
\begin{equation*}
\int_{0}^{\infty}k\hat{K}\left(k\right)\frac{d}{dk}\left(\left|\hat{f}_A\left(k\right)\right|^{2}+\left|\hat{f}_A\left(-k\right)\right|^{2}\right)dk=\lambda\int_{A}x\frac{d}{dx}\left|f\left(x\right)\right|^{2}dx,
\end{equation*}
which results in \eqref{eq:energy} upon integration by parts and another use of the Parseval's identity.

The derivation \eqref{eq:energy} given here is deduced by the author for a specific kernel in \cite[Lem 2.1.1]{Po16}, but, as it was later discovered, it is similar to the proof of a lemma formulated in \cite{GoPa75} for a concrete bandlimited kernel.

The second step of the proof of the simplicity of the spectrum is to show that the non-vanishing of boundary values is incompatible with an eigenvalue having more than one eigenfunction. This final step goes as in \cite{SlPo61, GoPa75}.
First of all, we note that if $f$ satisfies \eqref{eq:main_eq}, then, by parity of $K$, it follows that so do its even and odd parts: $f_e\left(x\right):=\frac{1}{2}\left(f\left(x\right)+f\left(-x\right)\right)$, $f_o\left(x\right):=\frac{1}{2}\left(f\left(x\right)-f\left(-x\right)\right)$.
Now suppose that $\phi$ and $\psi$ are linearly independent solutions of \eqref{eq:main_eq} corresponding to the same $\lambda$. Differentiating both sides of the equation for $\phi$, we integrate by parts the integral term to obtain
\[
\lambda \phi^{\prime}\left(x\right)=\int_{A}K\left(x-t\right)\phi^{\prime}\left(t\right)dt-\phi\left(a\right)K\left(x-a\right)+\phi\left(-a\right)K\left(x+a\right), \hspace{1em}x\in A.
\]
We now integrate this equation against $\psi$ to deduce that
\begin{equation}
\lambda\left[\phi\left(a\right)\psi\left(a\right)-\phi\left(-a\right)\psi\left(-a\right)\right]=0.\label{eq:bnd_restr}
\end{equation}
By taking $\phi_e$, $\psi_o$ instead of $\phi$, $\psi$ in \eqref{eq:bnd_restr}, we conclude from \eqref{eq:energy} that either $\phi_e\equiv 0$ or $\psi_o\equiv 0$.
Likewise, taking $\phi_o$, $\psi_e$ instead of $\phi$, $\psi$ in \eqref{eq:bnd_restr}, we deduce that $\phi_o\equiv 0$ or $\psi_e\equiv 0$. In other words, both $\phi$ and $\psi$ should be either odd or even, but if it is so, define $w\left(x\right):=\phi\left(x\right)\psi\left(a\right)-\phi\left(a\right)\psi\left(x\right)$. Clearly, we have $w\left(\pm a\right)=0$ and hence, invoking again  \eqref{eq:energy}, we conclude that $w\equiv 0$ contradicting the assumed linear independence of $\phi$ and $\psi$.

Finally, we note that simplicity of the spectrum combined with the parity and real-valuedness of $K$ implies that each eigenfunction of $\mathfrak{K}$ is either odd or even and can be chosen to be real-valued.
\end{proof}

In what follows, we will assume real-valuedness and the normalisation $\left\|f\right\|_{L^2\left(A\right)}=1$.

\section{Asymptotic solution for $\left|A\right|\gg 1$}
\label{sec:large_A}

We will obtain the main constructive result of this Section as a consequence of the following representation theorem. 

\begin{theorem}\label{thm:sol_repr}
Suppose that $K$ satisfies Assumptions \ref{ass:K_par}-\ref{ass:FT_monot_decay}. If $f$, $\lambda$ solve \eqref{eq:main_eq}, then they also satisfy for $x\in A$
\begin{equation}\label{eq:f_sol_repr}
f\left(x\right)=\frac{2\pi i \hat{K}\left(k_{0}\right)}{\hat{K}^{\prime}\left(k_{0}\right)}\hat{\phi}_{+}\left(k_{0}\right)\left(e^{2\pi ik_{0}\left(a-x\right)}\pm e^{2\pi ik_{0}\left(a+x\right)}\right)+\int_{\mathbb{R}}\hat{\phi}_{+}\left(k\right)\hat{T}_{0}\left(k\right)\left(e^{2\pi ik\left(a-x\right)}\pm e^{2\pi ik\left(a+x\right)}\right)dk,
\end{equation}
where the choice of sign in $\pm$ corresponds to even (plus) and odd (minus) parity of $f$, $k_0>0$ is a reparametrisation of $\lambda$ such that $\hat{K}\left(k_0\right)=\lambda$ and the following quantities are introduced
\begin{equation}
\hat{\phi}_{+}\left(k\right):=\frac{if\left(a\right)}{2\pi\left(k+i\right)X_{+}\left(k\right)}+q\left(k\right),\hspace{1em}q\left(k\right):=\frac{P_{+}\left[\hat{e}_{+}/X_{-}\right]\left(i\right)-P_{+}\left[\hat{e}_{+}/X_{-}\right]\left(k\right)}{4\pi^2\left(k^2+1\right)X_{+}\left(k\right)},\label{eq:phi_q_def}
\end{equation}

\begin{equation}
\hat{T}_{0}\left(k\right):=\frac{\hat{K}\left(k\right)}{\hat{K}\left(k\right)-\hat{K}\left(k_{0}\right)}-\frac{2k_{0}}{\left(k^{2}-k_{0}^{2}\right)\left(\log\hat{K}\right)^{\prime}\left(k_{0}\right)},\hspace{1em}
\hat{e}_+\left(k\right):=\int_{0}^{\infty}e^{2\pi ikx}e_a\left(x\right)dx,\label{eq:T0_def}
\end{equation}

\begin{equation}
e_a\left(x\right):=\int_{0}^{\infty}\left[N_{0}^{\prime\prime}\left(x+t+2a\right)+4\pi^2 k^2_{0}N_{0}\left(x+t+2a\right)-\frac{8\pi^3 k_0\left(k_0^2+1\right)}{\left(\log \hat{K}\right)^{\prime}\left(k_0\right)}e^{-2\pi\left(x+t+2a\right)}\right]f\left(t+a\right)dt,\label{eq:e_a_def}
\end{equation}
\begin{equation}
N_{0}\left(x\right):=\int_{\mathbb{R}}e^{-2\pi ikx}\left[\frac{\hat{K}\left(k\right)}{\hat{K}\left(k\right)-\hat{K}\left(k_0\right)}-\frac{2 k_{0}}{\left(\log \hat{K}\right)^{\prime}\left(k_0\right)}\left(\frac{1}{k^{2}-k_{0}^2}-\frac{1}{k^2+1}\right)\right]dk,\label{eq:N0_def}
\end{equation}
\begin{equation}
X_{\pm}\left(k\right):=G^{1/2}\left(k\right)\exp\left(\pm\frac{i}{2}\mathcal{H}\left[\log G\right]\left(k\right)\right),\hspace{1em}G\left(k\right):=\frac{\hat{K}\left(k_0\right)\left(k^{2}-k_0^{2}\right)}{\left(\hat{K}\left(k_0\right)-\hat{K}\left(k\right)\right)\left(k^{2}+1\right)}.\label{eq:X_pm_G_def}
\end{equation}
\end{theorem}
Note that operators $\mathcal{H}$, $P_+$, $P_-$ are as defined in the notational part of Section \ref{sec:asss_and_props}. 

\begin{proof}
	The proof is based on the reduction of \eqref{eq:main_eq} to an equation on $\mathbb{R}\backslash A$ whose solution admits an integral representation leading to \eqref{eq:f_sol_repr}.
	
	In view of continuity and the monotonicity of $\hat{K}\left(\left|k\right|\right)$ (Assumptions \ref{ass:K_reg}-\ref{ass:FT_monot_decay}) and the eigenvalue bounds from Lemma \ref{lem:props}, for a given $\lambda\in\left(0,\hat{K}\left(0\right)\right)$, there exists a unique $k_{0}>0$ such that $\hat{K}\left(\pm k_{0}\right)=\hat{K}\left(k_{0}\right)=\lambda$. 
	
	Let us define for $x\in\mathbb{R}\backslash A$
	\begin{equation}
	f\left(x\right)=\frac{1}{\lambda}\int_{A}K\left(x-t\right)f\left(t\right)dt, \label{eq:f_ext_from_eq}
	\end{equation}
	and note that, due to the regularity and decay assumptions on $K$, we immediately have $f\in C^2\left(\mathbb{R}\right)$ and $f^{\left(j\right)}\left(x\right)=\mathcal{O}\left(1/x^\alpha\right)$, $j\in\left\{0,1,2\right\}$, $\alpha>1$, for $x\gg 1$.
	
	With smooth extension \eqref{eq:f_ext_from_eq}, equation \eqref{eq:main_eq} becomes valid for all $x\in\mathbb{R}$,
	and we can convolve it with the kernel function
	$R_{\epsilon}\left(x\right)$ defined for each $\epsilon>0$ as
	\begin{equation}
	R_{\epsilon}\left(x\right):=\int_{\mathbb{R}}\frac{e^{-2\pi ikx}\hat{K}\left(k\right)}{\hat{K}\left(k\right)-\lambda-i\epsilon}dk,\hspace{1em}x\in\mathbb{R}.\label{eq:R_def_expl}
	\end{equation}
	We have $R_\epsilon\in C^{2}\left(\mathbb{R}\right)$ and $R_{\epsilon}^{\left(j\right)}\left(x\right)=\mathcal{O}\left(1/x^{\alpha}\right)$, $j\in\left\{0,1,2\right\}$, $\alpha>1$, for $\left|x\right|\gg1$ since for $\epsilon>0$ formula \eqref{eq:R_def_expl} is the inverse Fourier transform of a function whose regularity and decay properties are the same as $\hat{K}$. 
	
	It is easy to verify in the Fourier domain that the following identity holds
	\begin{equation}
	\int_{\mathbb{R}}R_{\epsilon}\left(x-\tau\right)K\left(\tau\right)d\tau-\lambda R_{\epsilon}\left(x\right)=K\left(x\right)+i\epsilon R_{\epsilon}\left(x\right),\hspace{1em}x\in\mathbb{R},\label{eq:R_def_impl}
	\end{equation}
	and hence convolution of both sides of \eqref{eq:main_eq} with $R_\epsilon$ gives
	\begin{equation}
	f\left(x\right)+i\dfrac{\epsilon}{\lambda}\int_{A}R_{\epsilon}\left(x-t\right)f\left(t\right)dt=\int_{\mathbb{R}\backslash A}R_{\epsilon}\left(x-t\right)f\left(t\right)dt,\hspace{1em}x\in\mathbb{R},\label{eq:f_prelim_line_eq}
	\end{equation}
	where we exchanged the integration order (employing the regularity and decay of $R_\epsilon$ discussed above) and reused \eqref{eq:main_eq}. 
	
	Let us rewrite \eqref{eq:R_def_expl} as
	\begin{align}
	R_{\epsilon}\left(x\right)=& \int_{\mathbb{R}}e^{-2\pi ikx}\Biggl[\frac{\hat{K}\left(k\right)}{\hat{K}\left(k\right)-\lambda-i\epsilon}-\frac{2\lambda\left(k_{0}+i\epsilon/\hat{K}^{\prime}\left(k_{0}\right)\right)}{\hat{K}^{\prime}\left(k_{0}\right)}\Biggl(\frac{1}{k^{2}-\left(k_{0}+i\epsilon/\hat{K}^{\prime}\left(k_{0}\right)\right)^{2}}\Biggr.\Biggr.\nonumber \\
	& \Biggl.\Biggl.-\frac{1}{k^2+1}\Biggr)\Biggr]dk-\frac{2\pi\lambda\left(k_{0}+i\epsilon/\hat{K}^{\prime}\left(k_{0}\right)\right)}{\hat{K}^{\prime}\left(k_{0}\right)}e^{-2\pi\left|x\right|}-\frac{2\pi i\lambda}{\hat{K}^{\prime}\left(k_{0}\right)}e^{-2\pi\left(ik_{0}-\frac{\epsilon}{\hat{K}^{\prime}\left(k_{0}\right)}\right)\left|x\right|}\nonumber\\
	=:& N_{\epsilon}\left(x\right)+M_{\epsilon}\left(x\right)+L_{\epsilon}\left(x\right).\label{eq:R_def_expl2}
	\end{align}

	The integrand in the $N_{\epsilon}$ term of \eqref{eq:R_def_expl2} is as smooth as $\hat{K}$ for any $\epsilon\geq0$ due to the cancellation of the singularities at $k=\pm k_{0}$
	appearing in the limit $\epsilon\searrow 0$, we thus have $N_{\epsilon}\left(x\right)=\mathcal{O}\left(1/x^{\alpha}\right)$, $\alpha>1$, for $\left|x\right|\gg1$. 
	Moreover, the decay of the integrand is also inherited by $\hat{K}$ though this decay cannot be faster than $\mathcal{O}\left(1/k^4\right)$ for $\left|k\right|\gg1$ due to the presence of the second term inside the square bracket. 
	
	We observe that the second term in the left-hand side of (\ref{eq:f_prelim_line_eq}) vanishes as $\epsilon\searrow0$. Indeed, the passage to the limit under the integral sign is justified by the dominated convergence theorem applicable since $R_\epsilon$ is continuous for all $\epsilon\geq 0$ and the integration variable $t$ is varying in a finite range. This would yield the solution for the original problem (for a suitable set of values of $\lambda$) once we know the smooth extension of this solution onto the half-line $x>a$ defined by \eqref{eq:f_ext_from_eq}. Precisely, using that each solution $f$ has a certain parity (Lemma \ref{lem:props}), we have
	\begin{align}
	f\left(x\right)&=\lim _{\epsilon\rightarrow 0^{+}} \int_{\mathbb{R}\backslash A}R_{\epsilon}\left(x-t\right)f\left(t\right)dt\nonumber \\
	&=\lim _{\epsilon\rightarrow 0^{+}}\int_{a}^{\infty}\left(R_{\epsilon}\left(x-t\right)\pm R_{\epsilon}\left(x+t\right)\right)f\left(t\right)dt,\hspace{1em}x\in A,\label{eq:f_prelim_sol}
	\end{align}
	where the plus and minus signs correspond to the even and odd sets of solutions, respectively.
	
	We now embark on finding this smooth half-line extension. To this end, we return to equation (\ref{eq:f_prelim_line_eq})
	and apply to its both sides the differential operator $D_{\epsilon}:=\frac{d}{dx}+2\pi\left(ik_{0}-\frac{\epsilon}{\hat{K}^{\prime}\left(k_{0}\right)}\right)\text{sgn }x$.
	Regularity and decay of the integrand allows differentiation under the integral sign. Note that since $\left|x\right|^\prime=\text{sgn }x$, we have for $x$, $t>a$,
	\begin{align*}
	D_{\epsilon}\left[R_{\epsilon}\left(x+t\right)\right]=&D_{\epsilon}\left[N_{\epsilon}\left(x+t\right)\right]+2\pi\left(ik_{0}-\frac{\epsilon}{\hat{K}^{\prime}\left(k_{0}\right)}\right)\left(\text{sgn }x-\text{sgn}\left(x+t\right)\right)L_{\epsilon}\left(x+t\right)\\
	&+2\pi\left[\left(ik_{0}-\frac{\epsilon}{\hat{K}^{\prime}\left(k_{0}\right)}\right)\text{sgn }x-\text{ sgn}\left(x+t\right)\right]M_{\epsilon}\left(x+t\right)\\
	=&N_{\epsilon}^{\prime}\left(x+t\right)+2\pi ik_{0}N_{\epsilon}\left(x+t\right)+2\pi\left(ik_{0}-\frac{\epsilon}{\hat{K}^{\prime}\left(k_{0}\right)}-1\right)M_{\epsilon}\left(x+t\right),
	\end{align*}
	\vspace{-15pt}
	\begin{align*}
	D_{\epsilon}\left[R_{\epsilon}\left(x-t\right)\right]=&D_{\epsilon}\left[N_{\epsilon}\left(x-t\right)\right]+2\pi\left(ik_{0}-\frac{\epsilon}{\hat{K}^{\prime}\left(k_{0}\right)}\right)\left(1-\text{sgn}\left(x-t\right)\right)L_{\epsilon}\left(x-t\right)\\
	&+2\pi\left(ik_{0}-\frac{\epsilon}{\hat{K}^{\prime}\left(k_{0}\right)}-\text{sgn}\left(x-t\right)\right)M_{\epsilon}\left(x-t\right).
	\end{align*}
	The resulting integrand is hence bounded uniformly in a positive neighborhood of $\epsilon=0$ whereas $f\left(t\right)$ is smooth and decays for large $t$ in an integrable manner (see \eqref{eq:f_ext_from_eq}). This allows us to pass to the limit as $\epsilon\searrow 0$. As before, this eliminates the term that involves values of $f$ on $A$ due to the uniform boundedness of the integrand $D_{\epsilon}\left[R_{\epsilon}\left(x-t\right)\right]$ for $t\in A$ and the factor $\epsilon$ in front of it. 
	Therefore, denoting $D_0$, $N_0$ the quantities $D_\epsilon$, $N_\epsilon$ evaluated at $\epsilon=0$, we obtain for $x>a$
	\begin{align}
	D_0\left[f\left(x\right)\right]=&\int_{a}^{\infty}\Biggl[D_{0}\left[N_{0}\left(x-t\right)\right]-\frac{4\pi^{2}k_{0}\lambda}{\hat{K}^{\prime}\left(k_{0}\right)}\left(ik_{0}-\text{sgn}\left(x-t\right)\right)e^{-2\pi\left|x-t\right|}\Biggr.\nonumber\\
	&+\Biggl.\frac{4\pi^{2}k_{0}\lambda}{\hat{K}^{\prime}\left(k_{0}\right)}\left(1-\text{sgn}\left(x-t\right)\right)e^{-2\pi ik_{0}\left|x-t\right|}\Biggr]f\left(t\right)dt\pm r_a\left(x\right),	\label{eq:f_prelim_eq}
	\end{align}
	where $r_a\left(x\right):=\int_{a}^{\infty}\left[N_{0}^{\prime}\left(x+t\right)+2\pi ik_{0}N_{0}\left(x+t\right)-\frac{4\pi^2 k_0\lambda\left(ik_0-1\right)}{\hat{K}^{\prime}\left(k_0\right)}e^{-2\pi\left(x+t\right)}\right]f\left(t\right)dt$.
	
	We observe that the obtained equation is very close to a classical half-line Wiener-Hopf
	equation with the exception that it contains a differential operator in the left-hand side and the term $r_a$ which is an integral transform with a sum kernel (known also as Hankel-type kernel or cross-correlation transform). Moreover, the kernel in the integral operator in \eqref{eq:f_prelim_eq} is non-symmetric and oscillatory at infinity. We will first deal with the latter issue.
	It is convenient to set $\phi\left(x\right):=f\left(x+a\right)$, $x>0$, and rewrite \eqref{eq:f_prelim_eq} as an equation on $x>0$
	\begin{align}
	\phi^{\prime}\left(x\right)+2\pi ik_{0}\phi\left(x\right)=&\int_{0}^{\infty}D_{0}\left[N_{0}\left(x-t\right)\right]\phi\left(t\right)dt-\frac{4\pi^{2}k_{0}\lambda}{\hat{K}^{\prime}\left(k_{0}\right)}\int_{0}^{x}\left(ik_0-1\right)e^{-2\pi\left(x-t\right)}\phi\left(t\right)dt\nonumber\\
	&+\frac{4\pi^{2}k_{0}\lambda}{\hat{K}^{\prime}\left(k_{0}\right)}\int_{x}^{\infty}\left(2e^{2\pi ik_{0}\left(x-t\right)}-\left(ik_0+1\right)e^{2\pi\left(x-t\right)}\right)\phi\left(t\right)dt\pm r_a\left(x+a\right).\label{eq:phi_prelim_eq}
	\end{align}
	In order to eliminate the oscillatory component from the kernel function, we apply another differential operator $\overline{D}_{0}=\frac{d}{dx}-2\pi ik_{0}$ to both sides of \eqref{eq:phi_prelim_eq}. After simplification, this yields for $x>0$
	\begin{equation}
	\phi^{\prime\prime}\left(x\right)+4\pi^{2}k_{0}^2\phi\left(x\right)=\int_{0}^{\infty}\left[\left(\frac{d^{2}}{dx^{2}}+4\pi^{2}k_{0}^{2}\right)N_{0}\left(x-t\right)-\frac{8\pi^3 k_0\lambda\left(k_0^2+1\right)}{\hat{K}^{\prime}\left(k_0\right)}e^{-2\pi\left|x-t\right|}\right]\phi\left(t\right)dt\pm e_a\left(x\right),\label{eq:phi_halfline_eq}
	\end{equation}
where $e_a\left(x\right):=\overline{D}_{0}r_a\left(x+a\right)$ which can be written more explicitly as in \eqref{eq:e_a_def}.
	
	We now extend the standard Wiener-Hopf approach \cite{No88,GaCh78} to deal with integro-differential equation \eqref{eq:phi_halfline_eq} treating 
	$e_a$ as a known inhomogeneous term. 
	
	Let us introduce a solution-dependent function 
	\[
	\psi\left(x\right):=-\chi_{\mathbb{R}_{-}}\left(x\right)\int_{0}^{\infty}\left[\left(\frac{d^{2}}{dx^{2}}+4\pi^{2}k_{0}^{2}\right)N_{0}\left(x-t\right)-\frac{8\pi^3 k_0\lambda\left(k_0^2+1\right)}{\hat{K}^{\prime}\left(k_0\right)}e^{-2\pi\left|x-t\right|}\right]\phi\left(t\right)dt.
	\]
	One can then easily verify the validity on $\mathbb{R}$ of the following extended
	version of \eqref{eq:phi_halfline_eq}
	\begin{align*}
	\chi_{\mathbb{R}_{+}}\left(x\right)\phi^{\prime\prime}\left(x\right)= & \int_{\mathbb{R}}\left(N_{0}^{\prime\prime}\left(x-t\right)+4\pi^{2}k_{0}^{2}N_{0}\left(x-t\right)-\frac{8\pi^3 k_0\lambda\left(k_0^2+1\right)}{\hat{K}^{\prime}\left(k_0\right)}e^{-2\pi\left|x-t\right|}\right)\chi_{\mathbb{R}_{+}}\left(t\right)\phi\left(t\right)dt\\
	&-4\pi^{2}k_{0}^2\chi_{\mathbb{R}_{+}}\left(x\right)\phi\left(x\right)+\psi\left(x\right)\pm\chi_{\mathbb{R}_{+}}\left(x\right)e_a\left(x\right).
	\end{align*}
	Denoting $\hat{\phi}_{+}\left(k\right):=\int_{0}^{\infty}e^{2\pi ikx}\phi\left(x\right)dx$,  $\hat{e}_{+}\left(k\right):=\int_{0}^{\infty}e^{2\pi ikx}e_a\left(x\right)dx$, 
	we take Fourier transform of both sides of the obtained equation and, after few cancellations and a rearrangement of the terms, we arrive at
	\begin{equation}
	\dfrac{4\pi^{2}\lambda\left(k^{2}-k_{0}^{2}\right)}{\lambda-\hat{K}\left(k\right)}\hat{\phi}_{+}\left(k\right)+\hat{\psi}\left(k\right)=-\phi^{\prime}\left(0\right)+2\pi ik\phi\left(0\right)\mp\hat{e}_{+}\left(k\right),\hspace{1em}k\in\mathbb{R},\label{eq:prelim_RH_pbm}
	\end{equation}
	where we used the convolution theorem and the fact that $\phi\left(x\right)$, $\phi^\prime\left(x\right)$ decay to zero as $x\rightarrow \infty$.
	
	By the Paley-Wiener theory, \eqref{eq:prelim_RH_pbm} can be viewed as a multiplicative Riemann-Hilbert
	problem for the unknown functions $\Phi_{+}\left(k\right):=4\pi^{2}\left(k^{2}+1\right)\hat{\phi}_{+}\left(k\right)$ and 
	$\Phi_{-}\left(k\right):=\psi\left(k\right)$ analytic in the upper
	and lower half-planes, respectively, and satisfying the following compatibility
	condition on their common boundary of analyticity
	\begin{equation}
	G\left(k\right)\Phi_{+}\left(k\right)+\Phi_{-}\left(k\right)=W\left(k\right)\mp\hat{e}_{+}\left(k\right),\hspace{1em}k\in\mathbb{R},\label{eq:RH_pbm}
	\end{equation}
	where $W\left(k\right):=-\phi^{\prime}\left(0\right)+2\pi ik\phi\left(0\right)$
	and 
	\begin{equation*}
	G\left(k\right):=\frac{\lambda\left(k^{2}-k_{0}^{2}\right)}{\left(\lambda-\hat{K}\left(k\right)\right)\left(k^{2}+1\right)}.
	\end{equation*}
	
	Note that the coefficient $G\left(k\right)$ is a piecewise differentiable (as is $\hat{K}$) positive
	non-vanishing function on the real line such that $G\left(k\right)\rightarrow1$ as
	$k\rightarrow\pm\infty$, hence $\log G\left(k\right)$ is single-valued
	and we can perform a standard procedure of finding its Wiener-Hopf
	factorisation $G\left(k\right)=X_{+}\left(k\right)X_{-}\left(k\right)$, $k\in \mathbb{R}$,
	where $X_{+}$, $X_{-}$ are analytic non-vanishing functions defined
	in the upper and the lower half-planes, respectively \cite{GaCh78,No88}.
	Such a factorisation can be achieved by the Plemelj-Sokhotskii formulas for boundary
	values of the Cauchy integral of $\log G\left(k\right)$, namely,
	\begin{align*}
	\log X_{\pm}\left(k\right)=P_{\pm}\left[\log G\right] \left(k\right)=\frac{1}{2}\left(\log G\left(k\right)\pm i\mathcal{H}\left[\log G\right]\left(k\right)\right)\nonumber\\
	\Rightarrow\hspace{1em}X_{\pm}\left(k\right)=G^{1/2}\left(k\right)\exp\left(\pm\dfrac{i}{2}\mathcal{H}\left[\log G\right]\left(k\right)\right),\hspace{1em}k\in\mathbb{R},
	\end{align*}
	where we used projectors $P_{+}$ and $P_{-}$ onto subspaces of functions defined on the real line which are traces of functions analytic in the upper and lower half-planes, respectively, and the Hilbert transform operator $\mathcal{H}$, all defined in the beginning of Section \ref{sec:asss_and_props}.
	
	The obtained factorisation $G\left(k\right)=X_{+}\left(k\right)X_{-}\left(k\right)$ together with the decomposition
	\[
	\hat{e}_{+}\left(k\right)/X_{-}\left(k\right)=P_{+}\left[\hat{e}_{+}/X_{-}\right]\left(k\right)+P_{-}\left[\hat{e}_{+}/X_{-}\right]\left(k\right),\hspace{1em}k\in\mathbb{R},
	\]
	reduce (\ref{eq:RH_pbm}) to the following relation 
	\[
	X_{+}\left(k\right)\Phi_{+}\left(k\right)+P_{+}\left[\hat{e}_{+}/X_{-}\right]\left(k\right)=\dfrac{W\left(k\right)-\Phi_{-}\left(k\right)}{X_{-}\left(k\right)}\mp P_{-}\left[\hat{e}_{+}/X_{-}\right]\left(k\right),\hspace{1em}k\in\mathbb{R}.
	\]
	Since this relation expresses equality of boundary values of analytic functions
	in the upper and lower half-planes along the common boundary of analyticity,
	we conclude that these are the restrictions of one entire function.
	By the generalized Liouville theorem, this function can only be a polynomial. We will now deduce that this polynomial is of degree $1$.
	Indeed, as can be seen from integration by parts in the Fourier integral, we have 
	\begin{equation}
	\hat{\phi}_{+}\left(k\right)=-\frac{1}{2\pi ik}\phi\left(0\right)+o\left(1/k\right),\hspace{1em}\left|k\right|\gg1,\label{eq:phi_hat_asympt}
	\end{equation}
	where smallness of the second term is furnished by the Riemann-Lebesgue lemma since $\phi^{\prime}\left(x\right)\in L^1\left(\mathbb{R}\right)$.
	Moreover, $\left|X_{\pm}\left(k\right)\right|=G^{1/2}\left(k\right)\rightarrow 1$
	as $k\rightarrow \pm\infty$. This leads us to conclude that
	\[
	X_{+}\left(k\right)\Phi_{+}\left(k\right)+P_{+}\left[\hat{e}_{+}/X_{-}\right]\left(k\right)=c_{0}+c_{1}k,\hspace{1em}\text{Im } k\geq 0,
	\]
	for some constants $c_{0}$, $c_{1}\in\mathbb{C}$, and hence, we obtain
	\begin{equation}
	\hat{\phi}_{+}\left(k\right)=\frac{c_{0}+c_{1}k-P_{+}\left[\hat{e}_{+}/X_{-}\right]\left(k\right)}{4\pi^{2}\left(k^{2}+1\right)X_{+}\left(k\right)},\hspace{1em}\text{Im } k\geq 0.\label{eq:phi_hat_sol_prelim}
	\end{equation}
	The denominator in \eqref{eq:phi_hat_sol_prelim} vanishes at $k=i$ whereas $\hat{\phi}_{+}\left(k\right)$
	is analytic in the upper half-plane, consequently, we must have that $c_{0}=-ic_{1}+P_{+}\left[\hat{e}_{+}/X_{-}\right]\left(i\right)$. Furthermore, by matching the leading order asymptotic behavior of \eqref{eq:phi_hat_sol_prelim} with \eqref{eq:phi_hat_asympt}, it follows that $c_{1}=2\pi i\phi\left(0\right)$. Altogether, this implies that \eqref{eq:phi_hat_sol_prelim} can be written as \eqref{eq:phi_q_def}.
	
	Recalling that $\phi\left(x\right)=f\left(x+a\right)$ for $x>0$ and $\lambda=\hat{K}\left(k_{0}\right)$, we return to \eqref{eq:f_prelim_sol},
	and use the Parseval's identity for inner products and the following rewrite of \eqref{eq:R_def_expl}
	\begin{align*}
	\hat{R}_{\epsilon}\left(k\right)=&\frac{\hat{K}\left(k\right)}{\hat{K}\left(k\right)-\lambda-i\epsilon}-\frac{2\lambda\left(k_{0}+i\epsilon/\hat{K}^{\prime}\left(k_{0}\right)\right)}{\hat{K}^{\prime}\left(k_{0}\right)}\frac{1}{k^{2}-\left(k_{0}+i\epsilon/\hat{K}^{\prime}\left(k_{0}\right)\right)^{2}}\\
	&+\frac{2\lambda}{\hat{K}^{\prime}\left(k_{0}\right)}\frac{k_{0}+i\epsilon/\hat{K}^{\prime}\left(k_{0}\right)}{\left(k-k_{0}-i\epsilon/\hat{K}^{\prime}\left(k_{0}\right)\right)\left(k+k_{0}+i\epsilon/\hat{K}^{\prime}\left(k_{0}\right)\right)}.
	\end{align*}
	This gives for $x\in A$
	\begin{align}
	f\left(x\right)=&\lim_{\epsilon\rightarrow0^{+}}\int_{\mathbb{R}}\left(e^{2\pi ik\left(a-x\right)}\pm e^{2\pi ik\left(a+x\right)}\right)\overline{\hat{R}_{\epsilon}\left(k\right)}\hat{\phi}_{+}\left(k\right)dk\label{eq:efct_sol_prelim}\\
	=&\int_{\mathbb{R}}\left(e^{2\pi ik\left(a-x\right)}\pm e^{2\pi ik\left(a+x\right)}\right)\left(\frac{\hat{K}\left(k\right)}{\hat{K}\left(k\right)-\lambda}-\frac{2\lambda}{\hat{K}^{\prime}\left(k\right)}\frac{k_{0}}{k^{2}-k_{0}^{2}}\right)\hat{\phi}_{+}\left(k\right)dk\nonumber\\
	&+\frac{2\lambda}{\hat{K}^{\prime}\left(k\right)}\lim_{\epsilon\rightarrow0^{+}}\int_{\mathbb{R}}\frac{\left(e^{2\pi ik\left(a-x\right)}\pm e^{2\pi ik\left(a+x\right)}\right)\left(k_{0}-i\epsilon/\hat{K}^{\prime}\left(k_{0}\right)\right)}{\left(k-k_{0}+i\epsilon/\hat{K}^{\prime}\left(k_{0}\right)\right)\left(k+k_{0}-i\epsilon/\hat{K}^{\prime}\left(k_{0}\right)\right)}\hat{\phi}_{+}\left(k\right)dk,\nonumber
	\end{align}
	where the passage to the limit $\epsilon\searrow 0$ in the first term is justified by the dominated convergence theorem (the integrand is nonsingular and absolutely integrable).
	Such a limit passage under the integral sign cannot be carried out for the second integral term, however, we observe that, due to the analyticity and the uniform decay of the integrand in the entire upper half-plane (i.e. for $x\in A$ the exponential factors are uniformly small for $\left|k\right|\gg1$ with any $\text{Im } k>0$), the integral can be evaluated explicitly using residue calculus, and it is equal to
	\[
	2\pi ik_{0}\frac{\left(e^{2\pi i\left(k_{0}-i\epsilon/\hat{K}^{\prime}\left(k_{0}\right)\right)\left(a-x\right)}\pm e^{2\pi i\left(k_{0}-i\epsilon/\hat{K}^{\prime}\left(k_{0}\right)\right)\left(a+x\right)}\right)}{2\left(k_{0}-i\epsilon/\hat{K}^{\prime}\left(k_{0}\right)\right)}\hat{\phi}_{+}\left(k_{0}-i\epsilon/\hat{K}^{\prime}\left(k_{0}\right)\right).
	\]
	
	Therefore, recalling that $\hat{K}\left(k_0\right)=\lambda$, and using the definition of $\hat{T}_0$ in  \eqref{eq:T0_def}, \eqref{eq:efct_sol_prelim} yields \eqref{eq:f_sol_repr}.	

\end{proof}

Let us explain why representation \eqref{eq:f_sol_repr}, despite its complicated form, appears to be useful.

First, suppose that we had $q\equiv 0$ in \eqref{eq:phi_q_def}. Then, \eqref{eq:f_sol_repr} would give the solution for \eqref{eq:main_eq} explicitly up to a constant $f\left(a\right)$ and yet unknown $k_0$. Evaluation of \eqref{eq:f_sol_repr} at $x=a$ would lead to elimination of the factor $f\left(a\right)\neq 0$ (recall Lemma \ref{lem:props}) yielding a transcendental equation to solve for $k_0>0$ 
\begin{align}
1+\frac{2\pi \hat{K}\left(k_{0}\right)\left(1\pm e^{4\pi ik_{0}a}\right)}{\hat{K}^{\prime}\left(k_{0}\right)\left(k_0+i\right)X_{+}\left(k_{0}\right)}=\frac{i}{2\pi}\int_{\mathbb{R}}\frac{\hat{T}_{0}\left(k\right)\left(1\pm e^{4\pi ika}\right)}{\left(k+i\right)X_{+}\left(k\right)}dk,\label{eq:char_eqs}
\end{align}
where the sign choice in $\pm$ would correspond to even or odd sets of eigenfunctions.
In view of the one-to-one correspondence between $k_0>0$ and $\lambda\in\left(0,\hat{K}\left(0\right)\right)$, \eqref{eq:char_eqs} is, in fact, a pair of characteristic equations for finding sets of admissible values $k_0$ and hence $\lambda$. Once an eigenvalue $\lambda=\hat{K}\left(k_0\right)$ is found, $k_0$ can be inserted back to \eqref{eq:f_sol_repr} to produce a corresponding eigenfunction.

Now we provide a motivation for the fact that even though it is not true that $q\equiv0$, this quantity is small when the interval $A$ is large. 

We start off by noticing that $N_0$ is continuous and decays at infinity as fast as $K$ since \eqref{eq:N0_def} is the Fourier integral of an absolutely integrable function which is as smooth as $\hat{K}$. The same is true for $N_0^{\prime\prime}$ which essentially differs from \eqref{eq:N0_def} only by the presence of an extra factor $k^2$ under the integral sign. Altogether, we thus have $N_0\left(x\right)$, $N_0^{\prime\prime}\left(x\right)=\mathcal{O}\left(1/{\left|x\right|}^\alpha\right)$, $\alpha>1$, for $\left|x\right|\gg 1$.
Using this decay rate and the fact that the last term in the integrand of \eqref{eq:e_a_def} decays faster (exponentially), we can estimate for sufficiently large $a>0$ 
\begin{equation}\label{eq:e_a_ptwise}
\left|e_a\left(x\right)\right|\leq\frac{C_\lambda}{\left(x+a\right)^\alpha} \left\|f\left(\cdot+a\right)\right\|_{L^1\left(\mathbb{R}_{+}\right)}\leq \frac{C_\lambda}{\lambda \left(x+a\right)^\alpha} \left\|K\right\|_{L^1\left(\mathbb{R}\right)} \left\|f\right\|_{L^1\left(A\right)},\hspace{1em} x>0,
\end{equation}
where we employed the Young's inequality for convolutions applied to \eqref{eq:f_ext_from_eq}. Note that $C_\lambda>0$ is a constant that depends on $\lambda$ (as does $k_0$). By the Cauchy-Schwarz inequality, we can bound $\left\|f\right\|_{L^1\left(A\right)}\leq \left(2a\right)^{1/2} \left\|f\right\|_{L^2\left(A\right)}$,
and hence, taking into account the chosen normalisation $\left\|f\right\|_{L^2\left(A\right)}=1$, we have for $a\gg 1$
\begin{equation}
\left\|e_a\right\|_{L^2\left(\mathbb{R}_{+}\right)}=C_{\lambda}\cdot\mathcal{O}\left(\frac{1}{a^{\alpha-1}}\right).\label{eq:e_a_L2}
\end{equation}
Since, by the Parseval's identity, we have  $\left\Vert \hat{e}_{+}\right\Vert _{L^{2}\left(\mathbb{R}\right)}=\left\Vert e_{a}\right\Vert _{L^{2}\left(\mathbb{R}_{+}\right)}$, we can estimate
\[
\left\Vert P_{+}\left[\hat{e}_{+}/X_{-}\right]\right\Vert _{L^{2}\left(\mathbb{R}\right)}\le\left\Vert \hat{e}_{+}/X_{-}\right\Vert _{L^{2}\left(\mathbb{R}\right)}\leq\left\Vert \frac{1}{X_{-}}\right\Vert _{L^{\infty}\left(\mathbb{R}\right)}\left\Vert \hat{e}_{+}\right\Vert _{L^{2}\left(\mathbb{R}\right)}=\left\Vert \frac{1}{X_{-}}\right\Vert _{L^{\infty}\left(\mathbb{R}\right)}\left\Vert e_{a}\right\Vert _{L^{2}\left(\mathbb{R}_{+}\right)},
\]
\begin{align*}
\left|P_{+}\left[\hat{e}_{+}/X_{-}\right]\left(i\right)\right|&=\frac{1}{2\pi}\left|\int_{\mathbb{R}}\frac{\hat{e}_{+}\left(\tau\right)/X_{-}\left(\tau\right)}{\tau-i}d\tau\right|\leq\frac{1}{2\pi}\left\Vert \hat{e}_{+}/X_{-}\right\Vert _{L^{2}\left(\mathbb{R}\right)}\left(\int_{\mathbb{R}}\frac{d\tau}{\tau^{2}+1}\right)^{1/2}\\
&\leq\frac{1}{2\sqrt{\pi}}\left\Vert \frac{1}{X_{-}}\right\Vert _{L^{\infty}\left(\mathbb{R}\right)}\left\Vert e_{a}\right\Vert _{L^{2}\left(\mathbb{R}_{+}\right)}.
\end{align*}
Therefore, recalling \eqref{eq:X_pm_G_def} and using \eqref{eq:e_a_L2}, we obtain
\begin{equation}
\left\Vert q\right\Vert _{L^{2}\left(\mathbb{R}\right)}\leq\frac{1}{4\pi^{2}}\left(1+\frac{1}{2\sqrt{\pi}}\right)\left\Vert \frac{1}{G}\right\Vert _{L^{\infty}\left(\mathbb{R}\right)}\left\Vert e_{a}\right\Vert _{L^{2}\left(\mathbb{R}_{+}\right)}=C_{\lambda}\cdot\mathcal{O}\left(\frac{1}{a^{\alpha-1}}\right).\label{eq:q_estim}
\end{equation}

Continuity and the smallness of $q$ in $L^{\infty}\left(\mathbb{R}\right)$ norm can also be shown. The corresponding estimates hinge on the following two facts about Fourier and Hilbert transforms. First, if $F\in L^{\infty}\left(\mathbb{R}\right)$  such that $F\left(x\right)=\mathcal{O}\left(1/\left|x\right|^{\alpha}\right)$, $\alpha>1$, for $\left|x\right|\gg 1$, then $\hat{F}\in C^{0,\beta}\left(\mathbb{R}\right)$ with some $\beta\in\left(0,1\right)$. Second, if $F\in C^{0,\beta}\left(\mathbb{R}\right)\cap L^2\left(\mathbb{R}\right)$, $\beta\in\left(0,1\right)$, then $\mathcal{H}\left[F\right]\in C^{0,\beta}\left(\mathbb{R}\right)\cap L^2\left(\mathbb{R}\right)$. In particular, \eqref{eq:e_a_ptwise} yields that $\hat{e}_+\in C^{0,\beta}\left(\mathbb{R}\right)\cap L^2\left(\mathbb{R}\right)$ with some $\beta>0$ whereas \eqref{eq:X_pm_G_def} with $G\in C^{0,1}\left(\mathbb{R}\right)$ implies that $1/X_{-}\in C^{0,\beta}\left(\mathbb{R}\right)$. Therefore, using the definition of $P_+$, we conclude that $P_{+}\left[\hat{e}_{+}/X_{-}\right]\in C^{0,\beta}\left(\mathbb{R}\right)$. 

Characteristic equations \eqref{eq:char_eqs} are much simpler than those obtained by
the author in \cite{BaLePo19} for a particular kernel of the
considered class. However, we can make another observation that provides even further simplification. We notice that away from the endpoints of the interval $A$, i.e. for $x\in\left(-a+a^\gamma, a-a^\gamma\right)$ with any $0<\gamma<1$, the integral term in (\ref{eq:f_sol_repr}) is small for $a\gg1$ due to rapid oscillations. Indeed, this integral can be viewed as the Fourier transform of an integrable function evaluated at a large argument. By the Riemann-Lebesgue lemma, such an integral must decay to zero for large arguments since  $\hat{\phi}_+$, $\hat{T}_0$ are continuous and their product is absolutely integrable.
Therefore, when $x\in\left(-a+a^\gamma, a-a^\gamma\right)$, the non-integral oscillatory
term in (\ref{eq:f_sol_repr}) dominates for large $a$. On the other hand, according to Lemma \ref{lem:props}, we know that all eigenfunctions are real-valued (by the assumed normalisation).
This imposes a restriction on the complex phase of the constant in
front of the oscillatory function which is either $\cos\left(2\pi k_{0}x\right)$
(in case of even eigenfunctions) or $\sin\left(2\pi k_{0}x\right)$
(in case of odd eigenfunctions). Employing \eqref{eq:X_pm_G_def}, \eqref{eq:phi_q_def} and \eqref{eq:q_estim},
we are thus able to deduce the following versions of approximate characteristic equations
to be solved for $k_{0}>0$ (and, consequently, for $\lambda=\hat{K}$$\left(k_{0}\right)$)
\begin{equation}
2\pi k_{0}a-\arctan\dfrac{1}{k_{0}}-\dfrac{1}{2}\mathcal{H}\left[\log G\right]\left(k_{0}\right)\simeq \pi m,\hspace{1em}m\in\mathbb{Z},\hspace{1em}\text{(even eigenfunctions)}\label{eq:char_eq_even_simpl}
\end{equation}
\begin{equation}
2\pi k_{0}a-\arctan\dfrac{1}{k_{0}}-\dfrac{1}{2}\mathcal{H}\left[\log G\right]\left(k_{0}\right)\simeq \pi\left(m+\dfrac{1}{2}\right),\hspace{1em}m\in\mathbb{Z}.\hspace{1em}\text{(odd eigenfunctions)}\label{eq:char_eq_odd_simpl}
\end{equation}
Note that larger eigenvalues $\lambda$ correspond to the solutions with smaller positive values of $k_{0}$. 

For each solution $k_{0}$ of transcendental equations \eqref{eq:char_eqs}
or (\ref{eq:char_eq_even_simpl})-(\ref{eq:char_eq_odd_simpl}), the
corresponding eigenfunction would then be given by \eqref{eq:f_sol_repr} (with the choice of an appropriate
sign $\pm$). Using \eqref{eq:phi_q_def} and \eqref{eq:q_estim}, we have explicitly the following approximations valid for $a\gg 1$.\\
Even eigenfunctions:
\begin{align}
f\left(x\right)\simeq -\dfrac{2e^{2\pi ik_{0}a}\hat{K}\left(k_{0}\right)}{\left(k_{0}+i\right)X_{+}\left(k_{0}\right)\hat{K}^{\prime}\left(k_{0}\right)}\cos\left(2\pi k_{0}x\right)+\frac{i}{\pi}\int_{\mathbb{R}}\dfrac{e^{2\pi ika}\hat{T}_0\left(k\right)}{\left(k+i\right)X_{+}\left(k\right)}\cos\left(2\pi kx\right)dk,\label{eq:f_sol_even}
\end{align}
Odd eigenfunctions:
\begin{align}
f\left(x\right)\simeq \dfrac{2 i e^{2\pi ik_{0}a}\hat{K}\left(k_{0}\right)}{\left(k_{0}+i\right)X_{+}\left(k_{0}\right)\hat{K}^{\prime}\left(k_{0}\right)}\sin\left(2\pi k_{0}x\right)+\frac{1}{\pi}\int_{\mathbb{R}}\dfrac{e^{2\pi ika}\hat{T}_0\left(k\right)}{\left(k+i\right)X_{+}\left(k\right)}\sin\left(2\pi kx\right)dk.\label{eq:f_sol_odd}
\end{align}
Note that the pre-trigonometric factors in the first term of each (\ref{eq:f_sol_even})-(\ref{eq:f_sol_odd})
could be further simplified due to the validity of characteristic equations (\ref{eq:char_eq_even_simpl})-(\ref{eq:char_eq_odd_simpl}).

Altogether, the reasoning above suggests the following constructive description of the solution of \eqref{eq:main_eq} in the large interval case.

	For $k\in\mathbb{R}$, $\kappa\in\mathbb{R_{+}}$, let us set 
	\[
	\mathcal{G}\left(k,\kappa\right):=\frac{\hat{K}\left(\kappa\right)\left(k^{2}-\kappa^{2}\right)}{\left(\hat{K}\left(\kappa\right)-\hat{K}\left(k\right)\right)\left(k^{2}+1\right)},\hspace{1em} \mathcal{X}_{+}\left(k,\kappa\right):=\mathcal{G}^{1/2}\left(k,\kappa\right)\exp\left(\dfrac{i}{2}\mathcal{H}\left[\log\mathcal{G}\left(\cdot,\kappa\right)\right]\left(k\right)\right).
	\]
	
Denote as $\left(\kappa_{n}^{\left(e\right)}\right)_{n=1}$,
$\left(\kappa_{n}^{\left(o\right)}\right)_{n=1}$ all positive
solutions of two sets of transcendental equations 
\begin{equation}
\label{eq:transc_eq_even}
2\pi\kappa a-\arctan\dfrac{1}{\kappa}-\dfrac{1}{2}\mathcal{H}\left[\log\mathcal{G}\left(\cdot,\kappa\right)\right]\left(\kappa\right)=\pi m,\hspace{1em}m\in\mathbb{Z},
\end{equation}
\begin{equation}
\label{eq:transc_eq_odd}
2\pi\kappa a-\arctan\dfrac{1}{\kappa}-\dfrac{1}{2}\mathcal{H}\left[\log\mathcal{G}\left(\cdot,\kappa\right)\right]\left(\kappa\right)=\pi\left(m+\frac{1}{2}\right),\hspace{1em}m\in\mathbb{Z},
\end{equation}
respectively, assuming that these solutions are sorted in ascending order,
i.e. $\kappa_{n-1}^{\left(e\right)}<\kappa_{n}^{\left(e\right)}$,
$\kappa_{n-1}^{\left(o\right)}<\kappa_{n}^{\left(o\right)}$, $n\in\mathbb{N}_{+}$.
Eigenvalues corresponding to even and odd eigenfunctions are then given by
\[
\lambda_{n}^{\left(e\right)}=\hat{K}\left(\kappa_{n}^{\left(e\right)}\right)+\delta_{n,a}^{\left(e\right)},\hspace{1em}\lambda_{n}^{\left(o\right)}=\hat{K}\left(\kappa_{n}^{\left(o\right)}\right)+\delta_{n,a}^{\left(o\right)},
\]
for some constants $\delta_{n,a}^{\left(e\right)}$, $\delta_{n,a}^{\left(o\right)}\in\mathbb{R}$
such that $\delta_{n,a}^{\left(o\right)}$, $\delta_{n,a}^{\left(e\right)}\rightarrow0$
as $a\rightarrow+\infty$ for any fixed $n\in\mathbb{N}_{+}$. Corresponding
sets of even and odd eigenfunctions are, respectively, furnished by
\begin{align}
f_{n}^{\left(e\right)}\left(x\right)= & \dfrac{2\left(-1\right)^{n-1}}{\left[-2\kappa_{n}^{\left(e\right)}\left(\log\hat{K}\right)^{\prime}\left(\kappa_{n}^{\left(e\right)}\right)\right]^{1/2}}\cos\left(2\pi\kappa_{n}^{\left(e\right)}x\right)\label{eq:e_fcts_even}\\
& +\frac{i}{\pi}\int_{\mathbb{R}}\dfrac{e^{2\pi ika}}{\left(k+i\right)\mathcal{X}_{+}\left(k,\kappa_{n}^{\left(e\right)}\right)}\Biggl[\dfrac{\hat{K}\left(k\right)}{\hat{K}\left(k\right)-\hat{K}\left(\kappa_{n}^{\left(e\right)}\right)}\Biggr.\nonumber\\
& -\Biggl.\dfrac{2\kappa_{n}^{\left(e\right)}}{\left(k^{2}-\left(\kappa_{n}^{\left(e\right)}\right)^{2}\right)\left(\log\hat{K}\right)^{\prime}\left(\kappa_{n}^{\left(e\right)}\right)}\Biggr]\cos\left(2\pi kx\right)dk+\mathcal{E}_{n,a}^{\left(e\right)}\left(x\right),\nonumber
\end{align}
\begin{align}
f_{n}^{\left(o\right)}\left(x\right)= & \dfrac{2\left(-1\right)^{n-1}}{\left[-2\kappa_{n}^{\left(o\right)}\left(\log\hat{K}\right)^{\prime}\left(\kappa_{n}^{\left(o\right)}\right)\right]^{1/2}}\sin\left(2\pi\kappa_{n}^{\left(o\right)}x\right)\label{eq:e_fcts_odd}\\
& +\frac{1}{\pi}\int_{\mathbb{R}}\dfrac{e^{2\pi ika}}{\left(k+i\right)\mathcal{X}_{+}\left(k,\kappa_{n}^{\left(o\right)}\right)}\Biggl[\dfrac{\hat{K}\left(k\right)}{\hat{K}\left(k\right)-\hat{K}\left(\kappa_{n}^{\left(o\right)}\right)}\Biggr.\nonumber\\
& -\Biggl.\dfrac{2\kappa_{n}^{\left(o\right)}}{\left(k^{2}-\left(\kappa_{n}^{\left(o\right)}\right)^{2}\right)\left(\log\hat{K}\right)^{\prime}\left(\kappa_{n}^{\left(o\right)}\right)}\Biggr]\sin\left(2\pi kx\right)dk+\mathcal{E}_{n,a}^{\left(o\right)}\left(x\right)\nonumber,
\end{align}
with some error terms $\mathcal{E}_{n,a}^{\left(e\right)}$, $\mathcal{E}_{n,a}^{\left(o\right)}$ such that $\mathcal{E}_{n,a}^{\left(e\right)}$, 
$\mathcal{E}_{n,a}^{\left(o\right)}\rightarrow 0$ as $a\rightarrow+\infty$ in an appropriate norm for any fixed $n\in\mathbb{N}_{+}$.

Note that each eigenfunction in \eqref{eq:e_fcts_even}-\eqref{eq:e_fcts_odd} is defined by the $L^2$-normalisation up to a unimodular constant which, by the assumed real-valuedness, should be either $1$ or $-1$.

\section{Asymptotic solution for $\left|A\right|\ll 1$}
\label{sec:small_A}
The situation with equations on asymptotically small intervals is generally less exciting, but still has some theoretical and practical interest \cite{KnKe91,GrMi01}. 
As before, we assume the kernel function $K\left(x\right)$ to be even, real-valued (Assumption \ref{ass:K_par}), sufficiently smooth (see below) and such that all eigenvalues are simple (i.e. Assumptions \ref{ass:pos_def}-\ref{ass:FT_monot_decay}). 

It is convenient to start by rescaling. Let us denote $\varphi\left(x\right):=f\left(ax\right)$,
$\eta=\lambda/a$ and rewrite equation (\ref{eq:main_eq}) as

\begin{equation}
\int_{-1}^{1}K\left(a\left(x-t\right)\right)\varphi\left(t\right)dt=\eta\varphi\left(x\right),\hspace{1em}x\in\left(-1,1\right).\label{eq:main_eq_resc}
\end{equation}

Since the interval $A$ is small, a smooth kernel function can be well approximated
by a linear combination of very few elementary basis functions such as monomials. In this case, assuming that $K\in C^{M+1}\left(-a,a\right)$, the kernel function approximation by the Taylor series expansion reads $K\left(ax\right)=\sum_{m=0}^{M}\frac{K^{\left(2m\right)}\left(0\right)}{\left(2m\right)!}\left(ax\right)^{2m}$, where only even powers of $x$ enter the expansion due to the even parity of the kernel function. Inserting such an approximation in \eqref{eq:main_eq_resc} and using binomial formula on each term $\left(x-t\right)^{2m}$, it is evident that the resulting integral operator is of a rank $2M+1$ and hence cannot have more than $M+1$ even and $M$ odd eigenfunctions.
Therefore, to reproduce rich enough structure of the original operator, $M$ has to be taken sufficiently large. 

The idea of our approach is still to take advantage of the fact that
only few degrees of freedom are needed for a decent approximation
of the kernel function when $a$ is small, and attempt to improve the approximation of the corresponding integral operator. We pursue this goal by selecting an approximation class for kernel function such that, on
one hand, using each function in that class results in the approximation by a non-degenerate (infinite rank) compact operator and, on the other hand, the corresponding integral equation can either be solved explicitly or in a way more efficient than the original equation.

In particular, we select a kernel approximant from the following functional family 
\begin{equation}
K_{C,b,c}\left(ax\right):=C\frac{\sin\left(abx\right)}{\sin\left(acx\right)},\hspace{1em}x\in\left(-1,1\right),\label{eq:K_bc}
\end{equation}
defined for parameters $b$, $c\in\mathbb{R}\cup i\mathbb{R}$, $\left|b\right|\neq\left|c\right|$, $C\in\mathbb{C}$ chosen such that $K_{C,b,c}\left(ax\right)$ is real-valued for $x\in\left(-1,1\right)$ and, in case of $c\in\mathbb{R}$, $\left|c\right|<\pi/a$.

Integral operators with kernels from the class (\ref{eq:K_bc})
are very special since, as it is shown in \cite{Mo62}, they commute with a self-adjoint second-order differential operator,
a property that is very rare for smooth kernels \cite{Gr83,Wr90}. Due to this commutation, the integral and differential operators possess
a common set of eigenfunctions, and since the spectrum is simple,
the original problem is equivalent to solving a boundary-value problem for an ordinary differential equation. Namely, the eigenfunctions are those solutions of 
\begin{equation}
-\left(\left(1-\frac{\sin^{2}\left(acx\right)}{\sin^{2}\left(ac\right)}\right)\varphi^{\prime}\left(x\right)\right)^{\prime}+a^2\left(b^{2}-c^{2}\right)\frac{\sin^{2}\left(acx\right)}{\sin^{2}\left(ac\right)}\varphi\left(x\right)=\mu\varphi\left(x\right),\hspace{1em}x\in\left(-1,1\right),\label{eq:ODE_general}
\end{equation}
that are regular on $\left[-1,1\right]$. The condition of finiteness of the solution at the endpoints $x=\pm1$
restricts the values of the eigenparameter $\mu\in\mathbb{R}$ to an infinite discrete set. 

Note that despite sharing the same set of eigenfunctions,
eigenvalues $\left(\mu_{n}\right)_{n=1}^{\infty}$ of differential
operator are different from eigenvalues of the integral operator with the kernel $K_{C,b,c}$ and even more from eigenvalues $\left(\lambda_{n}\right)_{n=1}^{\infty}$ of the original problem (\ref{eq:main_eq_resc}). The latter can be obtained upon substitution into the integral equation. One way of doing it is by plugging the solution $\varphi$ in \eqref{eq:main_eq_resc} and integrating both sides against $\varphi$ on $\left(-1,1\right)$. This results in a numerically reliable formula 
\begin{equation}
\lambda_{n}=\frac{a\int_{-1}^{1}\varphi_{n}\left(x\right)\int_{-1}^{1}K_{C,b,c}\left(a\left(x-t\right)\right)\varphi_{n}\left(t\right)dtdx}{\left\Vert \varphi_{n}\right\Vert _{L^{2}\left(-1,1\right)}^{2}},\hspace{1em}n\in\mathbb{N}_{+}.\label{eq:lmbd_num}
\end{equation}

Now we proceed by tuning the parameters $C$, $b$, $c$ to the kernel function $K\left(x\right)$.  
To be more precise, let us 
consider the following two situations. 
\begin{enumerate}
\item Suppose that $K\in C^3\left(-a,a\right)$. Then, as the approximant we can take for $x\in\left(-1,1\right)$
\begin{equation}
\label{eq:small_A_K_apprx_1}
\mathcal{K}_c\left(ax\right):=\frac{K\left(0\right)c}{\sqrt{c^2-3K^{\prime\prime}\left(0\right)/K\left(0\right)}}\frac{\sin\left(\sqrt{c^2-3K^{\prime\prime}\left(0\right)/K\left(0\right)}ax\right)}{\sin\left(acx\right)},
\end{equation}
where constant $c\in \mathbb{R}\cup i\mathbb{R}$ can be chosen arbitrary with the only restriction that $\left|c\right|<\pi/a$ when $c\in\mathbb{R}$. 
\item Suppose that $K\in C^5\left(-a,a\right)$. Then, a suitable choice of the approximant for $x\in\left(-1,1\right)$ is 
\begin{equation}
\label{eq:small_A_K_apprx_2}
\mathcal{K}\left(ax\right):=K\left(0\right)\sqrt{\frac{5K\left(0\right)K^{\left(4\right)}\left(0\right)-9\left[K^{\prime\prime}\left(0\right)\right]^2}{5K\left(0\right)K^{\left(4\right)}\left(0\right)-21\left[K^{\prime\prime}\left(0\right)\right]^2}}\frac{\sin\left(\sqrt{\frac{5K\left(0\right)K^{\left(4\right)}\left(0\right)-21\left[K^{\prime\prime}\left(0\right)\right]^2}{4 K\left(0\right)K^{\prime\prime}\left(0\right)}}ax\right)}{\sin\left(\sqrt{\frac{5K\left(0\right)K^{\left(4\right)}\left(0\right)-9\left[K^{\prime\prime}\left(0\right)\right]^2}{4 K\left(0\right)K^{\prime\prime}\left(0\right)}}ax\right)}.
\end{equation}
\end{enumerate}

Note that we do not require the quantities under the square root signs in \eqref{eq:small_A_K_apprx_1} and \eqref{eq:small_A_K_apprx_2} to be positive since the same expressions appear outside and inside the sine functions and hence still make the approximant to the kernel to be real-valued.

Equations \eqref{eq:small_A_K_apprx_1} and \eqref{eq:small_A_K_apprx_2} follow from the Taylor expansion of \eqref{eq:K_bc}
\[
K_{C,b,c}\left(ax\right)=C\frac{b}{c}\left(1+\frac{c^{2}-b^{2}}{6}a^2x^{2}+\frac{3b^{4}-10b^{2}c^{2}+7c^{4}}{360}a^4x^{4}+o\left(x^{4}\right)\right),
\]
when two or three non-vanishing terms are retained, respectively, and then matched to the expansion of $K\left(ax\right)$ about $x=0$ to the corresponding order. This ensures coincindence of $K\left(0\right)$, $K^{\prime\prime}\left(0\right)$, $K^{\left(4\right)}\left(0\right)$ with the respective quantities of the approximant. 

We further point out one way to obtain an explicit approximation to the solution of \eqref{eq:main_eq_resc} in terms of special functions. To this effect, we focus on \eqref{eq:small_A_K_apprx_1} and observe that in the limiting case $c\rightarrow 0$ one has $\mathcal{K}_c\left(x\right)=K\left(0\right)\frac{\sin\left(\sqrt{-3K^{\prime\prime}\left(0\right)/K\left(0\right)}x\right)}{\sqrt{-3K^{\prime\prime}\left(0\right)/K\left(0\right)}x}$ which is the celebrated sinc kernel appearing in many contexts relevant to bandlimited signal theory. The main advantage of such an approximation is that the corresponding differential equation will have polynomial coefficients. Indeed, in the limit $c\rightarrow 0$ \eqref{eq:ODE_general} becomes
\[
\left(\left(1-x^{2}\right)\varphi^{\prime}\left(x\right)\right)^{\prime}+\left(\mu+3\frac{K^{\prime\prime}\left(0\right)}{K\left(0\right)}a^{2}x^{2}\right)\varphi\left(x\right)=0,\hspace{1em}x\in\left(-1,1\right),
\]
which is a well-studied equation (see e.g. \cite{OsRoXi13}, \cite{SlPo61}) whose solutions are bounded on $\left[-1,1\right]$
only for special values $\mu_{n}=\chi_{n-1}\left(\sqrt{-3\frac{K^{\prime\prime}\left(0\right)}{K\left(0\right)}}a\right)$,
$n\in\mathbb{N}_{+}$, and are known as prolate spheroidal wave functions
\begin{equation}
\varphi_{n}\left(x\right)=S_{0\left(n-1\right)}\left(\sqrt{-3\frac{K^{\prime\prime}\left(0\right)}{K\left(0\right)}}a,x\right),\hspace{1em} n\in\mathbb{N}_{+},\label{eq:e_fcts_prolate}
\end{equation}
(we follow the notation of \cite{SlPo61}) which thus asymptotically coincide, up to a sign, with the solutions of \eqref{eq:main_eq_resc}. Furthermore, according to \cite{SlPo61}, eigenvalues of the integral operator with the sinc kernel function $\mathcal{K}\left(ax\right)$ can be expressed in terms of another set of special functions known as radial spheroidal wave functions. Namely, we have 
\begin{equation}
\lambda_n=2K\left(0\right)a\left[R_{0\left(n-1\right)}^{\left(1\right)}\left(\sqrt{-3\frac{K^{\prime\prime}\left(0\right)}{K\left(0\right)}}a,1\right)\right]^2,\hspace{1em} n\in\mathbb{N}_{+},
\label{eq:lmbd_num_alt}
\end{equation}
(following the notation of \cite{SlPo61}), and hence it is not even necessary to use \eqref{eq:lmbd_num}.

\section{Numerical illustrations}
\label{sec:num_ill}
We illustrate the obtained results on a concrete equation by taking a particular kernel function satisfying the assumptions of Section \ref{sec:asss_and_props}. We compare asymptotic solutions of Sections \ref{sec:large_A}-\ref{sec:small_A} with the numerical ones which are computed on a mesh fine enough so that they can be considered as the reference. More precisely, we take $K\left(x\right)=1/\left(1+x^2\right)^{3/2}$ which is the kernel function encountered in a simple model of geological prospecting \cite{Ki11}. Note that  $\hat{K}\left(k\right)=4\pi\left|k\right|\mathsf{K}_1\left(2\pi\left|k\right|\right)>0$, $\hat{K}^\prime\left(k\right)=-8\pi^2\left|k\right|\mathsf{K}_0\left(2\pi\left|k\right|\right)<0$ for $k\in\mathbb{R}\backslash\left\{0\right\}$ where $\mathsf{K}_0$, $\mathsf{K}_1$ are modified Bessel functions, and hence this kernel satisfies Assumptions \ref{ass:K_par}-\ref{ass:FT_monot_decay}. As a numerical method to generate reference solutions, we use the Nystr\"{o}m method (see e.g. \cite[Sect 4.7]{Ha95}) with the Gauss-Legendre quadrature. Namely, as in the previous section, we consider a rescaled problem \eqref{eq:main_eq_resc} from which the original solution of \eqref{eq:main_eq} can be restored as $f\left(x\right)=\varphi\left(x/a\right)$, $x\in A$, and, fixing $N=100$, we replace the integral by the quadrature rule yielding 
\begin{equation}
\sum_{j=1}^{N}\omega_{j}\tilde{K}\left(x-x_{j}\right)\varphi_{j}=\lambda/a\varphi\left(x\right),\hspace{1em}x\in\left[-1,1\right],\label{eq:half_discr_pbm}
\end{equation}
where $\omega_{j}:=2/\left[\left(1-x_{j}^{2}\right)\left[P_{N}^{\prime}\left(x_{j}\right)\right]^2\right]=2\left(1-x_{j}^{2}\right)/\left[N^2\left[P_{N-1}\left(x_{j}\right)\right]^2\right]$ are the quadrature weights, $P_{N-1}\left(x\right)$ is the $\left(N-1\right)$-th Legendre polynomial, $x_j$ is the $j$-th root of $P_{N}\left(x\right)$, and we denoted $\varphi_{j}:=\varphi\left(x_j\right)$, $\tilde{K}\left(x\right):=K\left(ax\right)=1/\left(1+a^2x^2\right)^{3/2}$.

Evaluation of both sides of \eqref{eq:half_discr_pbm} at $x=x_l$, $l=1,\dots,N$, leads to a vector eigenvalue problem
\begin{equation}
\sum_{j=1}^N a\omega_{j}\tilde{K}\left(x_{l}-x_{j}\right)\varphi_j=\lambda\varphi_l,\hspace{1em}l=1,\dots,N,\label{eq:discr_pbm}
\end{equation}
whose solutions we denote as $\lambda_l$ and $\left(\varphi_j\right)_l$, $l,j=1,\dots,N$.

Once a discrete eigenfunction $\left(\varphi_j\right)_l$, $j=1,\dots,N$, and a corresponding eigenvalue $\lambda_l$ are found, they can be inserted in \eqref{eq:half_discr_pbm} to yield a continuous numerical eigenfunction $\varphi_l\left(x\right)$ as a smooth interpolant on $\left[-1,1\right]$. Normalisation constant is chosen such that $\left\Vert\varphi_l\right\Vert_{L^2\left(-1,1\right)}=1$.  

\begin{figure}[!h]
\centering\includegraphics[width=4.3in]{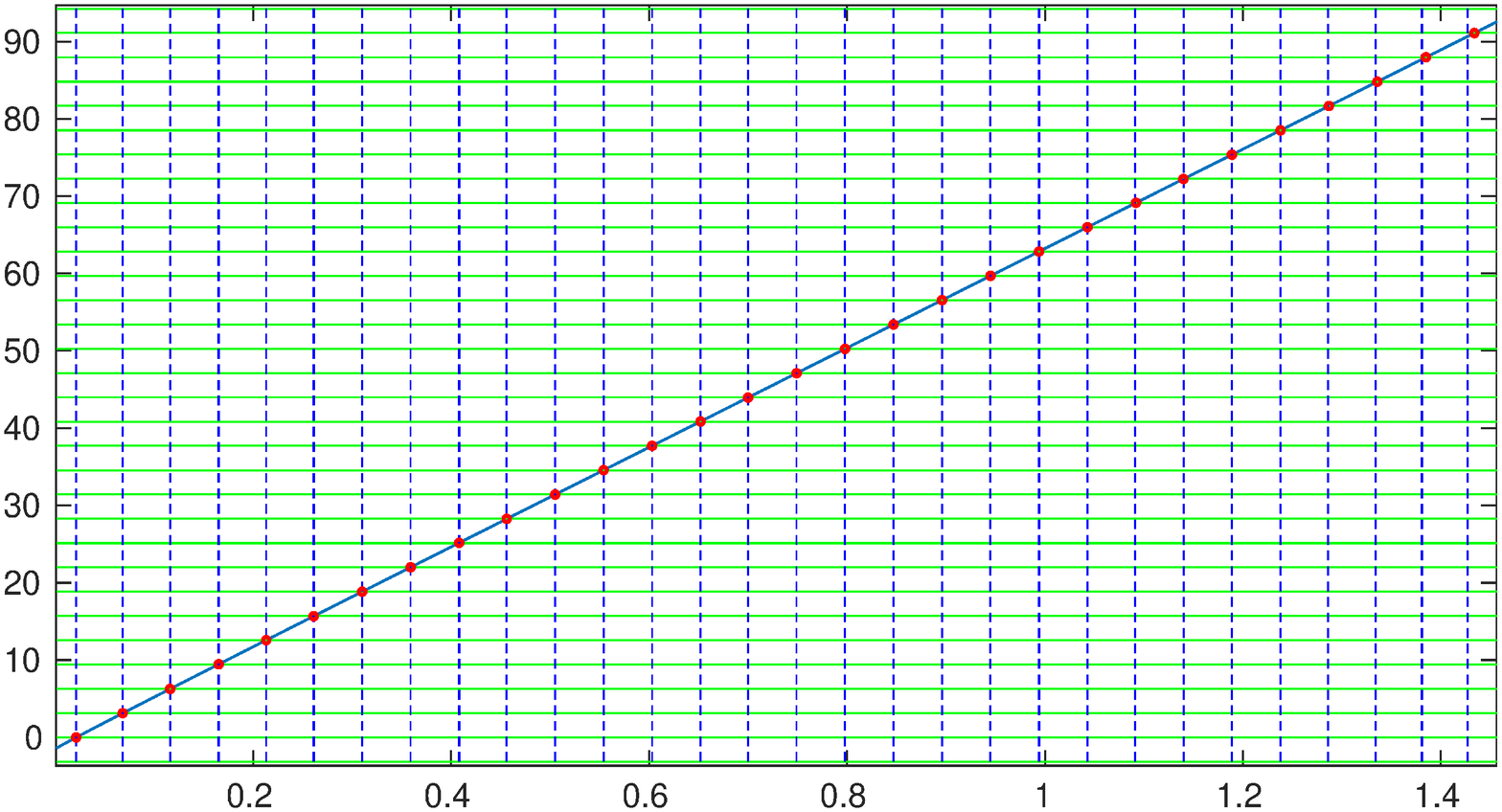}
\vskip6pt\vspace{10px}
\includegraphics[width=4.3in]{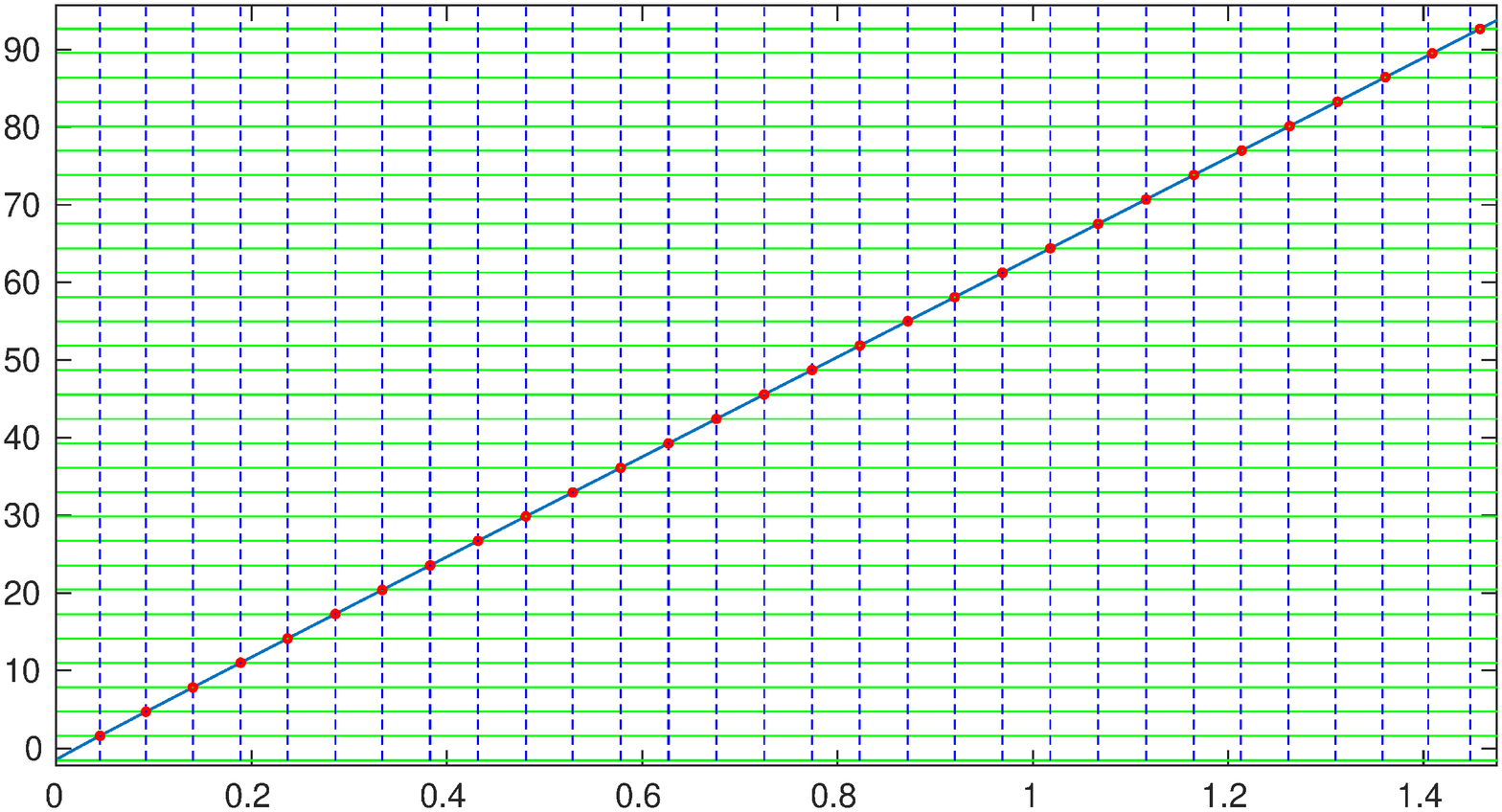}
\caption{Solving characteristic equations for even (top) and odd (bottom) part of the spectrum. Blue solid curve represents the left-hand side of \eqref{eq:char_eq_even_simpl} and \eqref{eq:char_eq_odd_simpl}, green lines illustrate a variety of right-hand side terms for $m=0$, $1$, $\dots$. Red dots are the intersection points detected by the algorithm, their abscissas thus correspond to the solutions of the characteristic equations and hence the asymptotically found eigenvalues. Blue vertical dash lines correspond to the reference eigenvalues (found numerically).}
\label{fig:char_eqs}
\end{figure}

We first verify the asymptotic results of Section \ref{sec:large_A} by taking $a=10$.
We illustrate graphically the solution of a set of characteristic equations \eqref{eq:transc_eq_even} for the even part of the spectrum on the top plot of Figure \ref{fig:char_eqs}.
We observe that for $m<0$ \eqref{eq:transc_eq_even} does not have a (positive) solution whereas for each fixed $m\in\mathbb{N}_{0}$ the corresponding equation has exactly one solution $\kappa^{\left(e\right)}_n$, $n:=m+1$, this is where the plotted in blue color left-hand side attains the value $\pi m$. This solution is marked by a red dot, and the resulting eigenvalue is $\lambda^{\left(e\right)}_n=\hat{K}\left(\kappa^{\left(e\right)}_n\right)$. Similarly, with the odd part of the spectrum where eigenvalues $\lambda^{\left(o\right)}_n$, $n=m+1$, are entailed by solutions of \eqref{eq:transc_eq_odd}.
Reference values of $\kappa^{\left(e\right)}_n$, $\kappa^{\left(o\right)}_n$ plotted as vertical lines are in unique correspondence to the values $\lambda^{\left(e\right)}_n=\lambda_{2n-1}$ and $\lambda^{\left(o\right)}_n=\lambda_{2n}$, $n=1,\dots,\left\lfloor N/2\right\rfloor+1$, which are found from the solution of \eqref{eq:discr_pbm}. In making this comparison, we used that even and odd eigenvalues interlace as they decrease in magnitude starting from the largest eigenvalue $\lambda_1$ corresponding to the even eigenfunction. The figure shows that red dots start deviating visibly from vertical lines only for larger values of $k$, i.e. for smaller eigenvalues, namely, those of the ordinal number $n\simeq 30$. Table \ref{tbl:a10} with relative eigenvalue and eigenfunction approximation errors confirms that the nearly first $60$ eigenvalues and eigenfunctions are well approximated by the derived large interval asymptotics \eqref{eq:transc_eq_even}-\eqref{eq:transc_eq_odd}, \eqref{eq:e_fcts_even}-\eqref{eq:e_fcts_odd} for $a=10$. 

For the case $a=0.1$, the results of the naive approximation strategy described in Section \ref{sec:small_A} are given in Table \ref{tbl:a01}. We take advantage of the existing code for computing prolate spheroidal functions \eqref{eq:e_fcts_prolate} provided in \cite{ZhJ96} which is converted to the MATLAB format using software f2matlab \cite{f2matlab}. Eigenvalues are computed from the prolate spheroidal functions using \eqref{eq:lmbd_num} rather than relying on the computation of radial spheroidal functions with \eqref{eq:lmbd_num_alt}, the latter approach is found to be less accurate numerically for the eigenvalues of very small magnitude. Evidently, the results are not as great as for the large interval case: we can obtain a reliable approximation only for the very first few eigenfunctions and eigenvalues. Note that in this case eigenvalues are very small, a fact that makes numerical computations challenging. On the other hand, for practical problems the contribution of higher eigenfunctions with small eigenvalues in the expansion of the measured data is expected to be negligible. It is noteworthy that the computation of eigenvalues breaks down before the computation of eigenfunctions. Since eigenvalues are computed from the knowledge of eigenfunctions, this drawback could be circumvented by using higher precision arithmetic or better evaluation strategies than \eqref{eq:lmbd_num}, but we do not pursue this any further here.

When comparison between the asymptotic and reference solutions is made, we recall that eigenfunctions are defined by \eqref{eq:e_fcts_even}-\eqref{eq:e_fcts_odd} and \eqref{eq:e_fcts_prolate}
only up to a sign (see a remark at the end of Section \ref{sec:large_A}). Therefore, we use the fact that $f_n\left(a\right)\neq 0$ and perform multiplication by $-1$ if necessary such that the signs of $f_n\left(a\right)$ for a reference eigenfunction and its asymptotic approximation match.   

\begin{figure}[!h]
\centering\includegraphics[width=4.3in]{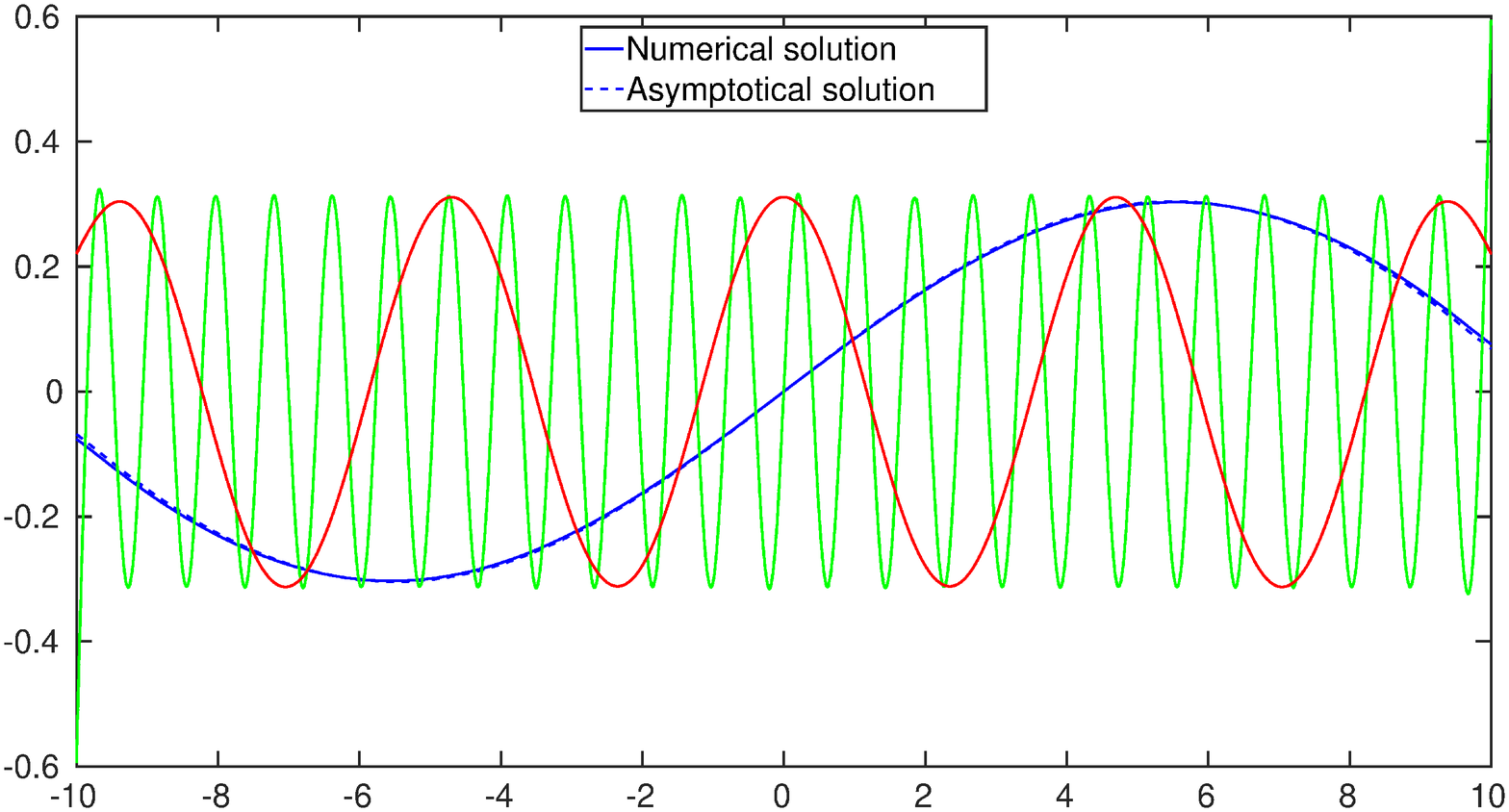}
\vskip6pt\vspace{10px}
\includegraphics[width=4.3in]{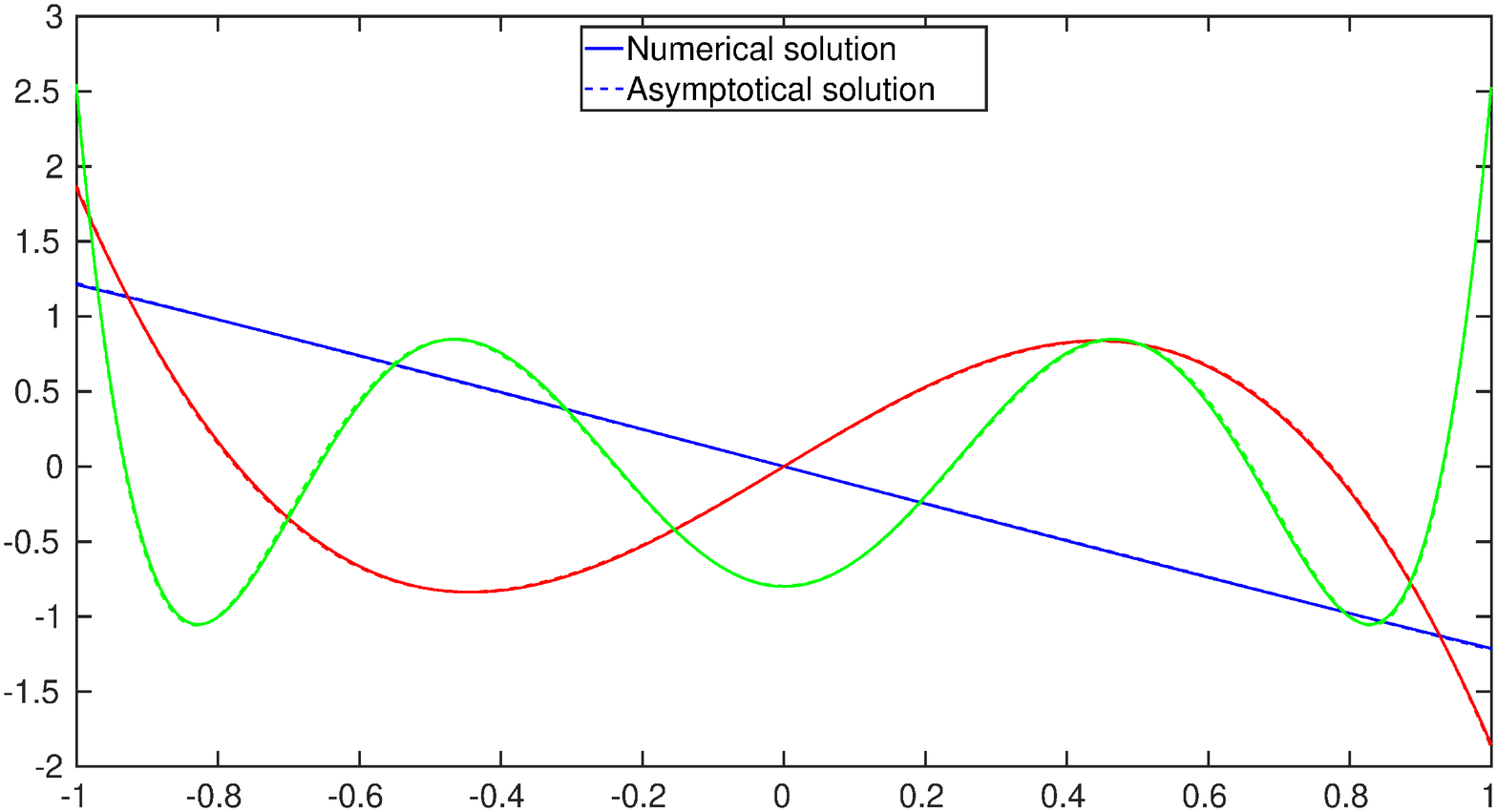}
\caption{Some eigenfunctions for $a=10$ (top): 1st odd (blue), 5th even (red), 25th odd (green), and for $a=0.1$ (bottom): 1st odd (blue), 2nd odd (red), 4th even (green). Dash lines correspond to asymptotic solutions, solid lines represent numerical reference solutions (visually indistinguishable).}
\label{fig:e_fcts}
\end{figure}

We show several eigenfunctions on Figure \ref{fig:e_fcts} by plotting both asymptotic and reference solutions, though they essentially cannot be distinguished visually. We see that in the case $a=10$, eigenfunctions deviate very little from trigonometric functions with this deviation being more pronounced near the endpoints only for the eigenfunctions of large number $n$. For the case $a=0.1$, the eigenfunctions resemble some orthogonal polynomials which is not surprising due to similarity of the differential equation for prolate spheroidal harmonics to that for Legendre polynomials.

\begin{figure}[!h]
\centering\includegraphics[width=4.3in]{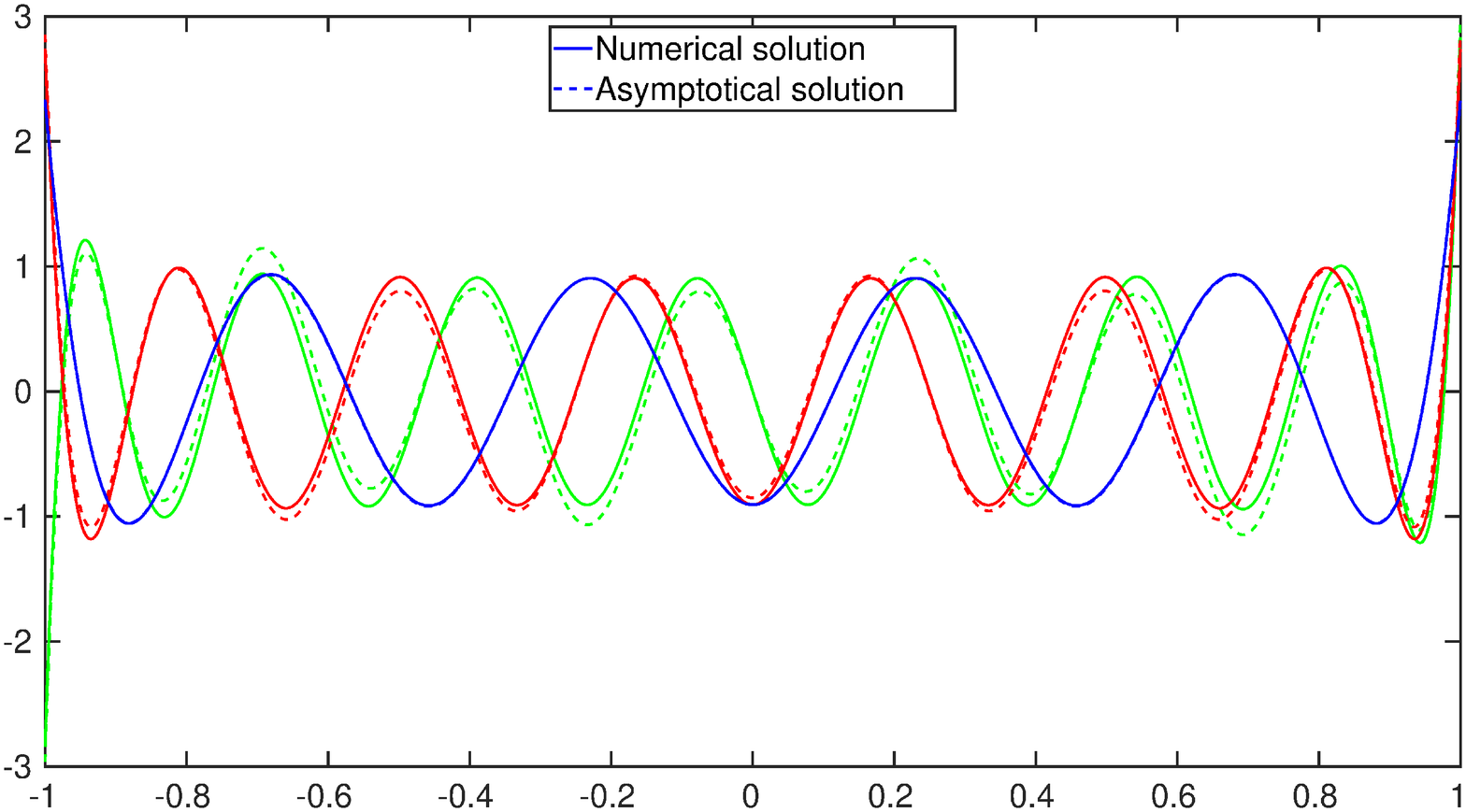}
\vskip6pt\vspace{10px}
\includegraphics[width=4.3in]{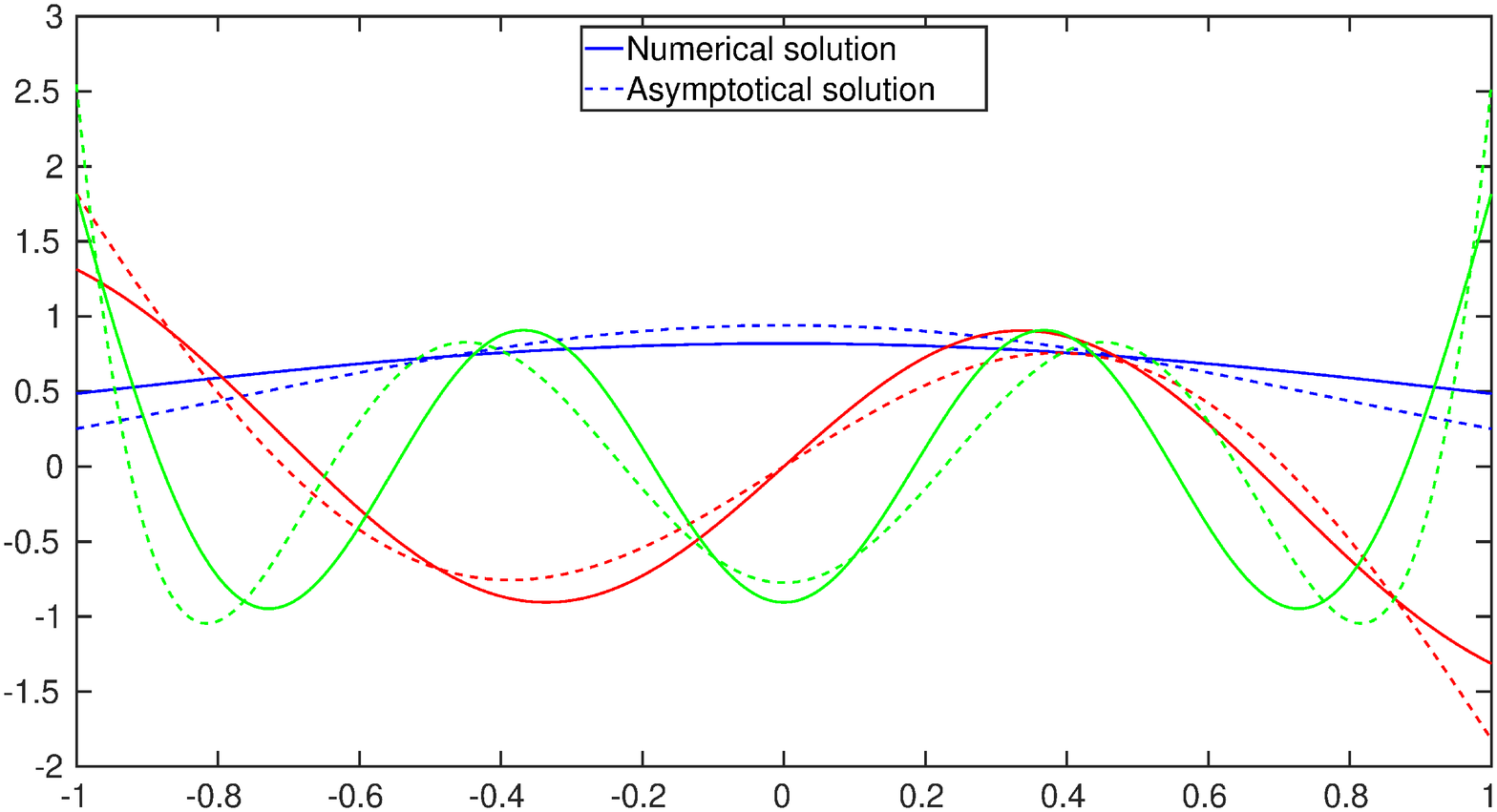}
\caption{Breakdown of asymptotic solutions for $a=1$. Eigenfunctions for large interval approximation (top): 6th even (blue), 7th even (red), 8th odd (green), and small interval approximation (bottom): 1st even (blue), 2nd odd (red), 3rd even (green). Dash lines correspond to asymptotic solutions, solid lines represent numerical reference solutions.}
\label{fig:asympt_breakdown}
\end{figure}

To illustrate the breakdown of the approximation to the eigenfunctions, we apply both asymptotic approaches described in Sections \ref{sec:large_A}-\ref{sec:small_A} to the case $a=1$ where none of them should be justifiably valid. We see on Figure \ref{fig:asympt_breakdown} (bottom) that even though the small interval asymptotic approximations to eigenfunctions deviate from the reference solutions immediately, they still reproduce qualitative behaviour well enough. On the other hand, it is remarkable that the large interval asymptotic approximations to eigenfunctions begin deviating visually from the reference solutions only starting from the number $n\simeq 14$. This deviation is shown on Figure \ref{fig:asympt_breakdown} (top).   

\begin{table}[!h]
\caption{Accuracy of eigenvalue and eigenfunction approximation (large interval asymptotics): $a=10$}
\label{tbl:a10}
\begin{tabular}{|l|l|l|l|l|l|l|}
\hline
$n$ & $\lambda_{n}^{\left(e\right)}$ & $\bigl|\delta\lambda_{n}^{\left(e\right)}\bigr|/\lambda_{n}^{\left(e\right)}$ & $\bigl\Vert \delta f_{n}^{\left(e\right)}\bigr\Vert _{L^{2}\left(A\right)}$ & $\lambda_{n}^{\left(o\right)}$ & $\bigl|\delta\lambda_{n}^{\left(o\right)}\bigr|/\lambda_{n}^{\left(o\right)}$ & $\bigl\Vert \delta f_{n}^{\left(o\right)}\bigr\Vert _{L^{2}\left(A\right)}$\\
\hline
1 & 1.95e+00 & 8.44e-05 & 1.60e-03 & 1.85e+00 & 1.17e-03 & 1.33e-02 \\\hline
2 & 1.72e+00 & 9.99e-04 & 9.44e-03 & 1.58e+00 & 8.70e-04 & 7.47e-03 \\\hline
3 & 1.44e+00 & 1.01e-03 & 8.12e-03 & 1.31e+00 & 1.50e-03 & 1.16e-02 \\\hline
4 & 1.17e+00 & 2.19e-03 & 1.64e-02 & 1.05e+00 & 1.35e-03 & 9.87e-03 \\\hline
5 & 9.41e-01 & 3.06e-04 & 2.19e-03 & 8.39e-01 & 9.58e-04 & 6.71e-03 \\\hline
6 & 7.46e-01 & 2.38e-03 & 1.64e-02 & 6.59e-01 & 1.06e-03 & 7.27e-03 \\\hline
7 & 5.84e-01 & 5.56e-04 & 3.75e-03 & 5.16e-01 & 2.31e-03 & 1.55e-02 \\\hline
8 & 4.53e-01 & 1.06e-03 & 7.06e-03 & 4.00e-01 & 8.37e-04 & 5.53e-03 \\\hline
9 & 3.50e-01 & 2.50e-03 & 1.64e-02 & 3.08e-01 & 4.84e-04 & 3.17e-03 \\\hline
10 & 2.71e-01 & 1.60e-03 & 1.04e-02 & 2.36e-01 & 1.68e-03 & 1.09e-02 \\\hline
11 & 2.07e-01 & 5.01e-04 & 3.24e-03 & 1.82e-01 & 2.74e-03 & 1.76e-02 \\\hline
12 & 1.58e-01 & 5.08e-04 & 3.29e-03 & 1.39e-01 & 1.80e-03 & 1.15e-02 \\\hline
13 & 1.21e-01 & 1.44e-03 & 9.18e-03 & 1.05e-01 & 9.26e-04 & 5.95e-03 \\\hline
14 & 9.16e-02 & 2.30e-03 & 1.46e-02 & 8.00e-02 & 1.08e-04 & 1.16e-03 \\\hline
15 & 6.98e-02 & 2.55e-03 & 1.62e-02 & 6.06e-02 & 6.63e-04 & 4.26e-03 \\\hline
16 & 5.29e-02 & 1.82e-03 & 1.16e-02 & 4.58e-02 & 1.40e-03 & 8.76e-03 \\\hline
17 & 3.99e-02 & 1.11e-03 & 7.31e-03 & 3.46e-02 & 2.10e-03 & 1.31e-02 \\\hline
18 & 3.01e-02 & 4.32e-04 & 3.53e-03 & 2.61e-02 & 2.79e-03 & 1.71e-02 \\\hline
19 & 2.27e-02 & 2.35e-04 & 2.43e-03 & 1.98e-02 & 2.34e-03 & 1.54e-02 \\\hline
20 & 1.71e-02 & 9.02e-04 & 5.38e-03 & 1.49e-02 & 1.67e-03 & 1.20e-02 \\\hline
21 & 1.28e-02 & 1.59e-03 & 8.77e-03 & 1.12e-02 & 9.74e-04 & 8.80e-03 \\\hline
22 & 9.65e-03 & 2.32e-03 & 1.20e-02 & 8.38e-03 & 1.99e-04 & 6.79e-03 \\\hline
23 & 7.24e-03 & 3.17e-03 & 1.49e-02 & 6.29e-03 & 7.24e-04 & 6.12e-03 \\\hline
24 & 5.46e-03 & 1.66e-03 & 2.06e-02 & 4.71e-03 & 1.93e-03 & 7.48e-03 \\\hline
25 & 4.09e-03 & 2.46e-04 & 2.10e-02 & 3.53e-03 & 3.64e-03 & 1.11e-02 \\\hline
26 & 3.06e-03 & 1.87e-03 & 2.45e-02 & 2.64e-03 & 6.29e-03 & 1.84e-02 \\\hline
27 & 2.29e-03 & 5.25e-03 & 3.43e-02 & 1.98e-03 & 1.07e-02 & 3.32e-02 \\\hline
28 & 1.71e-03 & 1.10e-02 & 5.53e-02 & 1.48e-03 & 1.82e-02 & 6.30e-02 \\\hline
29 & 1.28e-03 & 2.11e-02 & 9.73e-02 & 1.11e-03 & 2.64e-02 & 1.47e-01 \\\hline
30 & 9.56e-04 & 3.93e-02 & 1.77e-01 & 8.29e-04 & 5.69e-02 & 3.07e-01 \\\hline
\end{tabular}
\end{table}

\begin{table}[!h]
\caption{Accuracy of eigenvalue and eigenfunction approximation (large interval asymptotics): $a=1$}
\label{tbl:a1}
\begin{tabular}{|l|l|l|l|l|l|l|}
\hline
$n$ & $\lambda_{n}^{\left(e\right)}$ & $\bigl|\delta\lambda_{n}^{\left(e\right)}\bigr|/\lambda_{n}^{\left(e\right)}$ & $\bigl\Vert \delta f_{n}^{\left(e\right)}\bigr\Vert _{L^{2}\left(A\right)}$ & $\lambda_{n}^{\left(o\right)}$ & $\bigl|\delta\lambda_{n}^{\left(o\right)}\bigr|/\lambda_{n}^{\left(o\right)}$ & $\bigl\Vert \delta f_{n}^{\left(o\right)}\bigr\Vert _{L^{2}\left(A\right)}$\\
\hline
1 & 1.25e+00 & 5.08e-03 & 2.22e-03 & 4.95e-01 & 4.50e-03 & 1.37e-03 \\\hline
2 & 1.69e-01 & 1.92e-03 & 9.66e-04 & 5.44e-02 & 3.19e-03 & 1.76e-03 \\\hline
3 & 1.68e-02 & 1.27e-03 & 6.91e-04 & 5.09e-03 & 1.24e-03 & 7.96e-04 \\\hline
4 & 1.51e-03 & 3.38e-03 & 2.09e-03 & 4.48e-04 & 1.16e-03 & 6.39e-04 \\\hline
5 & 1.31e-04 & 1.67e-03 & 1.96e-03 & 3.80e-05 & 2.58e-03 & 3.83e-03 \\\hline
6 & 1.11e-05 & 1.35e-03 & 6.76e-03 & 3.18e-06 & 5.91e-04 & 1.24e-02 \\\hline
7 & 9.13e-07 & 1.18e-03 & 2.69e-02 & 2.61e-07 & 1.93e-03 & 4.96e-02 \\\hline
8 & 7.50e-08 & 2.61e-03 & 1.01e-01 & 2.13e-08 & 1.85e-03 & 1.98e-01 \\\hline
9 & 6.07e-09 & 3.80e-04 & 3.45e-01 & 1.72e-09 & 4.12e-04 & 6.66e-01 \\\hline
10 & 4.88e-10 & 1.81e-03 & 8.64e-01 & 1.39e-10 & 3.53e-03 & 1.18e+00 \\\hline
\end{tabular}
\end{table}

\begin{table}[!h]
\caption{Accuracy of eigenvalue and eigenfunction approximation (small interval asymptotics): $a=0.1$}
\label{tbl:a01}
\begin{tabular}{|l|l|l|l||l|l|l|l|}
\hline
$n$ & $\lambda_{n}$ & $\bigl|\delta\lambda_{n}\bigr|/\lambda_{n}$ & $\bigl\Vert \delta \varphi_{n}\bigr\Vert _{L^{2}\left(-1,1\right)}$ & $n$ & $\lambda_{n}$ & $\bigl|\delta\lambda_{n}\bigr|/\lambda_{n}$ & $\bigl\Vert \delta \varphi_{n}\bigr\Vert _{L^{2}\left(-1,1\right)}$\\
\hline
1 & 1.98e+00 & 9.82e-09 & 9.91e-05 & 2 & 1.90e-03 & 1.63e-05 & 4.00e-03 \\\hline
3 & 1.53e-05 & 8.19e-05 & 6.70e-03 & 4 & 1.39e-07 & 2.94e-01 & 1.00e-02 \\\hline
5 & 1.43e-09 & 9.42e-01 & 1.32e-02 & 6 & 1.46e-11 & 1.99e+00 & 1.63e-02 \\\hline
7 & 1.48e-13 & 3.58e+00 & 1.94e-02 & 8 & 1.48e-15 & 6.08e+00 & 2.73e-02 \\\hline
9 & 1.36e-17 & 4.15e-01 & 1.42e+00 & 10 & -2.91e-18 & 1.53e+00 & 1.43e+00 \\\hline
\end{tabular}
\end{table}

\section{Conclusion and outlook}
\label{sec:concl}
We have addressed the problem of asymptotic solution of a convolution integral equation on a finite interval under fairly general assumptions on the kernel function satisfied in many practical applications. We focussed on the homogeneous equation \eqref{eq:main_eq} and proposed two approaches to find its solutions, i.e. eigenfunctions and eigenvalues of the corresponding integral operator. The large interval approach is different from asymptotic methods previously available in the literature in that it does not require a fast decay of the kernel function at infinity. The proposed small interval approach furnishes a simple alternative to the classical Taylor expansion of the kernel function. Numerical results illustrate efficiency of both approaches even if the size of the interval is only moderately large or small. The number of well-approximated eigenfunctions and eigenvalues in the large (resp. small) interval approach grow with increasing (resp. decreasing) the size of the interval. In practical applications, kernels are often smooth functions which result in a fast (essentially exponential) decay of eigenvalues and hence small contributions of higher eigenfunctions to the expansion of the data. Therefore, in such contexts only first few eigenfunctions might be needed. In the large interval case, the decay rate is slower due to the presence of a factor which is small for large $a$ (see e.g. \cite[Thm 2]{Wi64}). However, as illustrated in Table \ref{tbl:a10}, for the interval of size $2a=20$ first $57$ eigenfunctions could be asymptotically found with $L^2\left(A\right)$ error smaller than $10\%$. To compare (see Table \ref{tbl:a1}), for the interval with $a=1$ we can get $15$ eigenfunctions with the same accuracy.

While the obtained results are remarkable, there are certainly things to improve and explore deeper.

The provided justification of the large interval approach is not complete. While we gave a proof of the representation theorem which is an essence of our method, the proof of the entailed approximation results was merely sketched.
In particular, the dependence of $\lambda$ (and equivalently $k_0$) on $a$ should be taken into account in establishing error bounds \eqref{eq:e_a_ptwise}-\eqref{eq:e_a_L2} and \eqref{eq:q_estim}.

The accuracy of the large interval approach could be improved to the exponential one (according to \eqref{eq:e_a_def}) if compactly supported or exponentially decaying extension of the kernel function $K$ outside the interval $\left(-2a,2a\right)$ is used when computing $\hat{K}$ and the dependent quantities $\mathcal{G}$, $\mathcal{X}_{+}$ entering the final solution formulae \eqref{eq:transc_eq_even}-\eqref{eq:transc_eq_odd}, \eqref{eq:e_fcts_even}-\eqref{eq:e_fcts_odd}. Of course, it is reasonable to argue that after such an extension, the methods for rapidly decaying kernels become applicable. However, the main analytical disadvantage of the extension procedure is that the Fourier transform of the resulting piecewise defined function is generally not expected to be computable in explicit form whereas such a computation often could be done (e.g. employing residue calculus) if the kernel is initially given by the restriction to $\left(-2a,2a\right)$ of a smooth function defined on the whole $\mathbb{R}$.   

The large interval approach was presented as a consequence of a solution representation theorem obtained due to a possibility of solving effectively an auxiliary integral equation on the complement of the interval. It might be feasible to iterate between solution of the problem on a complement of the interval and on the interval $A$ itself to produce the solution of arbitrary high accuracy for a fixed value of $a$. Convergence and practical utility of such an iterative scheme is yet to be investigated, but it is expected at least when the kernel (or its extension) decays exponentially at infinity (cf. \cite{Ki18}). 

Even though the obtained results allow solving \eqref{eq:main_eq_nonhom} through solution of \eqref{eq:main_eq} (as mentioned in Section \ref{sec:intro}), the large interval method could also be applicable directly to a inhomogeneous equation \eqref{eq:main_eq_nonhom} under less restrictive assumptions, at least when $\lambda>\lambda_1$ and $g$ is either even or odd function (see e.g. \cite{Gr64}).

Advantages of the solution of homogeneous equation \eqref{eq:main_eq} lie beyond solving \eqref{eq:main_eq_nonhom}. Indeed, the knowledge of eigenfunctions furnishes a spectral resolution of the integral operator. This spectral resolution is useful in obtaining a regularised solution of Fredholm integral equations of the first kind as well as operator equations involving functions of the "restriction-convolution-restriction" operator (e.g. see \cite{LePo17} when the square of such an operator appears). 

Due to the simplicity of the plots on Figures \ref{fig:char_eqs} and nearly periodic distribution of the solutions of characteristic equations, further investigation of a possibility of obtaining more explicit results is needed.   

It is noteworthy that numerical experiments with few kernel functions show that $\eqref{eq:char_eq_even_simpl}$ with $m=0$ yields the characteristic equation for the first even eigenvalue regardless of the size of the interval $A$. This might give a basis for further theoretical work on the estimates of the first eigenvalue (i.e. the spectral radius of an operator).

Since it is clear that convolution integral equations are continuous analogs of Toeplitz matrices, it is natural to see how the obtained results transplant into statements about the spectrum of some classes of large Toeplitz matrices (cf. \cite{BoBo15}).

In the small interval case, other members of the proposed approximation family \eqref{eq:K_bc} could be investigated. In particular, it would be curious to see which other ODEs from the family \eqref{eq:ODE_general} admit solutions in the form of conveniently computable special functions. Also, on the level of the integral equation approximation, an exponentially accurate or perhaps, with some additional effort, even an explicit form of the solution could be obtained for the kernel function $\text{sech }x\equiv 1/\cosh x$. 

\vskip6pt
\enlargethispage{20pt}



\section*{Acknowledgements}
The author would like to thank the Isaac Newton Institute for Mathematical Sciences, Cambridge, UK, for its support and hospitality during the programme "Bringing pure and applied analysis together via the Wiener-Hopf technique, its generalisations and applications" supported by EPSRC (EP/R014604/1), and personally Prof. Gennady Mishuris (Aberystwyth University) for the invitation and funding the author's participation in that programme. The author is also thankful to Prof. Sergey Rogosin (Belarusian State University) for useful comments on the work during the programme. Finally, the author would like to acknowledge the FWF-project I3538-N32 for the current employment and Sci-Hub service for providing the access to numerous papers cited in this work.



\vskip2pc


\end{document}